\documentclass[12pt, english]{article}
\usepackage[utf8]{inputenc}
\usepackage[12pt]{extsizes}
\usepackage[english]{babel}
\usepackage{enumitem}
\usepackage{dynkin-diagrams}
\usepackage{verbatim}
\usepackage{blindtext}
\usepackage{geometry}
\usepackage{amsmath}
\usepackage{amsthm}
\usepackage{amssymb}
\usepackage{dsfont}
\usepackage{mathrsfs}
\usepackage{graphicx}

\DeclareGraphicsExtensions{.pdf,.png,.jpeg,.jpg,.JPG}

\usepackage{fancyhdr} 
\usepackage{lastpage} 
\usepackage{extramarks} 
\usepackage{lipsum} 

\newtheorem{theorem}{Theorem}
\newtheorem{lemma}{Lemma}
\newtheorem{propos}{Proposition}
\newtheorem{definition}{Definition}
\newtheorem{remark}{Remark}

\newenvironment{Proof}
{\par\noindent{\bf Proof.}}
{\hfill$\scriptstyle\blacksquare$}

\begin{document}
\begin{center}
\textbf{Short $SL_2$-structures on simple Lie algebras}\\
\vspace{\baselineskip}
Roman O. Stasenko\\
\emph{Moscow Center of Fundamental and Applied Mathematics, Lomonosov Moscow State University, Moscow, Russia}\\
theromestasenko@yandex.ru\\
\vspace{\baselineskip}
\end{center}
\begin{footnotesize}
{\bf Abstract.}Throughout the papers of E.B. Vinberg some non-abelian gradings of simple Lie algebras were introduced and investigated, namely short $SO_3 -$ and $SL_3$ - structures. We study another kind of them -- short $SL_2$ - structures. The main results relate to one-to-one corresponding between such structures and some special Jordan algebras.
\end{footnotesize}
\newline
\newline
{\bf Key words.} Structural Lie algebras, representations of lie algebras, izotypic decomposition.
\newline
\newline
{\bf Introduction.} In the theory of Lie algebras, the classical Tits-Cantor-Koeher construction is well known. It allows using a simple Jordan algebra $J$ to construct a simple Lie algebra $\mathfrak{g}$, which has the following form:
\begin{equation}\label{ttk}
  \mathfrak{g} = \mathfrak{der}(J)\oplus\mathfrak{sl}_2(J).
\end{equation}

In this construction, the commutator of elements of the second term is defined as an ordinary commutator of matrices plus a certain addition from the first term. More precisely, if we identify $\mathfrak{sl}_2(J)$ with $\mathfrak{sl}_2(\mathbb{C})\otimes J$, then the commutator of the elements of the second term is calculated by the following formula:
$$[X\otimes a, Y\otimes b] = 2(X, Y)[L_a, L_b] + [X, Y]\otimes ab,\qquad\forall X, Y\in\mathfrak{sl}_2(\mathbb{C}), a, b\in J,$$
where $L_c$ denotes linear multiplication operator by the element $c\in J$, and

$$(X, Y) = \operatorname{tr}(XY),\quad\forall X, Y\in \mathfrak{sl}_2(\mathbb{C}).$$

The Tits-Cantor-Koeher construction can be interpreted as a linear representation of the group $SL_2\left(\mathbb{C}\right)$ by automorphisms of the algebra $\mathfrak{g}$, which decomposes into irreducible representations of dimensions 1 and 3. Equation (\ref{ttk}) is nothing but isotypic decomposition of the given representation.

A natural generalization of this linear representation is the following definition:
\begin{definition}
Let $S$ be a reductive algebraic group. The $S$-structure on the Lie algebra $\mathfrak{g}$ is called the homomorphism $\Phi:S\rightarrow\operatorname{Aut}(\mathfrak{g})$.
\end{definition}
In this terminology, the Tits-Cantor-Koeher construction will be called a very short $SL_2$-structure. A complete analysis of this construction can be found in \cite[section~3]{Vinberg_1}.

Some special cases of $S$-structures were studied in \cite{Vinberg_1} and \cite{Vinberg_2}. We will consider another case of $S$-structures, namely --- short $SL_2$-structures.

It is most natural to consider short $SL_2$-structures on semisimple Lie algebras, due to the remarkable algebraic properties of the latter. However, it is obvious that considering a short $SL_2$-structure on a semisimple Lie algebra reduces to considering such structures on every simple component of this algebra. Therefore, we will consider simple Lie algebra $\mathfrak{g}$ and $S = SL_2(\mathbb{C})$ everywhere.

The differential of the map $\Phi$ sets the linear representation of the Lie algebra $\mathfrak{sl}_2$ by differentiations of the Lie algebra $\mathfrak{g}$:
$$ \operatorname{d}\Phi: \mathfrak{sl}_2\rightarrow\mathfrak{der(g)}.$$

Since the algebra $\mathfrak{g}$ is simple, then $\mathfrak{der(g)}\simeq\mathfrak{inn(g)}\simeq\mathfrak{g}$. Hence, the image of the algebra $\mathfrak{sl}_2$ under the action of the map $\operatorname{d}\Phi$ in $\mathfrak{g}$ consists of adjoint operators of $\mathfrak{g}$.

The isotypic decomposition of this representation has the form:
$$\mathfrak{g} = \mathfrak{g}_0\oplus \overset{l}{\underset{i=1}{\bigoplus}}V_i\otimes J_i.$$
In this formula, $V_i$ denotes the space of $(i+1)$-dimensional irreducible representation of the algebra $\mathfrak{sl}_2$, and $J_i$ is a vector space on which the algebra $\mathfrak{sl}_2$ acts trivially, and which is responsible for the multiplicity of the occurrence of the irreducible representation $V_i$ in $\mathfrak{g}$ (namely, the dimension of $J_i$ is equal to this multiplicity). The subspace $\mathfrak{g}_0$ is an isotypic component corresponding to the one-dimensional representation of the algebra $\mathfrak{sl}_2$.
\begin{definition}
$SL_2$-the structure is called short if the representation is $\operatorname{d}\Phi$ decomposes into irreducible representations of dimensions 1, 2 and 3.
\end{definition}
An irreducible representation of an algebra $\mathfrak{sl}_2$ of dimension 2 is its tautological representation, and an irreducible representation of dimension 3 is the adjoint representation. Thus, for short $SL_2$--structures, the isotypic decomposition will have the form:
\begin{equation}\label{isot}
\mathfrak{g} = \mathfrak{g}_0\oplus\mathfrak{g}_1\oplus \mathfrak{g}_2,\quad\mathfrak{g}_1 = \mathbb{C}^2\otimes J_1,\quad\mathfrak{g}_2 = \mathfrak{sl}_2\otimes J_2.
\end{equation}
It is important to note that the cases of $J_1=0$ (which corresponds to a very short $SL_2$-structure) and $J_1\neq0$ differ significantly from each other. Next, we will assume everywhere that $J_1\neq0$.

Using the commutation operation on the Lie algebra $\mathfrak{g}$, as well as invariant scalar multiplication, the vector space $J_1$ is endowed with the structure of a symplectic space, and the space $J_2$ is endowed with the structure of the Jordan algebra of symmetric operators in $J_1$. Also, the algebra $\mathfrak{g}_0$ is a Lie subalgebra of the algebra $\mathfrak{sp}(J_1)$.

It turns out that among the short $SL_2$-structures with the given $J_1$ there is a maximum in the sense that the space $J_2$ and the Lie algebra $\mathfrak{g}_0$ in it are the maximum possible by inclusion. Namely, for this short $SL_2$-structure $\mathfrak{g}_0=\mathfrak{sp}(J_1)$, and the Jordan algebra of all symmetric operators of the symplectic space $J_1$ acts as the space $J_2$. We will call this Jordan algebra canonical. We will describe the maximum short $SL_2$-structure in the section $\ref{maxsl2str}$ of this article.

In the section $\ref{simpliejord}$ we will construct a correspondence between simple Lie algebras with a short $SL_2$-structure and so-called simple Lie-Jordan symplectic structures of the form $\left(J_1; J_2;\mathfrak{g}_0; \delta_0\right)$, where $J_2$ is a simple Jordan subalgebra of the algebra of all symmetric operators of the space $J_1$, $\mathfrak{g}_0$ is a reductive subalgebra in the algebra of all symplectic operators on the space $J_1$, and $\delta_0$ is some symmetric bilinear map, which we will describe below. This is the main result of the work formulated in the theorem $\ref{finth}$ of the subsection $\ref{maintheor}$.

In the section $\ref{classif}$, a complete classification of short $SL_2$-structures on simple Lie algebras will be carried out, indicating a simple symplectic Lie-Jordan structure that corresponds to each $SL_2$-structure.

The author expresses great gratitude for the formulation of the problem and scientific guidance to E. B. Vinberg, without whose advice and recommendations this article would not have seen the light, D. A. Timashev, whose clear instructions on editing allowed this text to take a slender and complete form, as well as D.I. Panyushev, who gave an accurate reference to one of the results of B. Kostant, required for the text of this article.

\section{Izotypic decomposition}\label{isdec}

In this section and everywhere else, all algebras, algebraic groups and linear representations are considered over the field $\mathbb{C}$.

\subsection{Commutational formulas}\label{comform}
 Consider an arbitrary reductive algebraic group $S$. For its irreducible representation, $\rho$ denote by $V_{\rho}$ the space of this representation. The following lemmas can be found in \cite[section 1,~subsection~1.2]{Vinberg_1}:

\begin{lemma}\label{lem1}
Let $\rho, \sigma$ and $\tau$ be irreducible representations of the group $S$. Assume that $\rho\otimes\sigma$ contains $\tau$ with multiplicity one and let
$$p: V_{\rho}\times V_{\sigma}\rightarrow V_{\tau}$$
be a fixed nonzero $S$-equivariant bilinear map (defined up to a scalar multiplication). Let $U_{\rho}, U_{\sigma}$ and $U_{\tau}$ be some vector spaces, on which the group $S$ acts trivially. Then every $S$-equivariant bilinear map
$$P: (V_{\rho}\times U_{\rho})\times (V_{\sigma}\otimes U_{\sigma}) \rightarrow V_{\tau}\otimes U_{\tau}$$
is given by the formula
$$P(a\otimes x, b\otimes y) = p(a,b)\otimes\nu(x,y),$$
where
$$\nu: U_{\rho}\times U_{\sigma}\rightarrow U_{\tau}$$
is some bilinear map.
\end{lemma}
\begin{lemma}\label{lem2}
Let $\rho$ and $\tau$ be irreducible representations of the group $S$. Assume that both $\operatorname{Sym}^2\rho$ and $\Lambda^2\rho$ contain $\tau$ with multiplicity at most one, and let
$$p: V_{\rho}\times V_{\rho} \rightarrow V_{\tau}\text{~and ~}q: V_{\rho}\times V_{\rho}\rightarrow V_{\tau}$$
be a fixed nonzero $S$-equivariant symmetric and skew-symmetric bilinear maps respectively (defined up to a scalar multiplication) in cases where $\operatorname{Sym}^2\rho$ or $\Lambda^2\rho$ contain $\tau$; otherwise, we put $p = 0\text{~or ~}q = 0)$. Let $U_{\rho}$ and $U_{\tau}$ be vector spaces, on which the group $S$ acts trivially. Then every $S$-equivariant skew-symmetric bilinear map
$$P: (V_{\rho}\otimes U_{\rho})\times (V_{\rho}\otimes U_{\rho}) \rightarrow V_{\tau}\otimes U_{\tau}$$
is given by the formula
$$P(a\otimes x, b\otimes y) = p(a, b)\otimes\phi(x, y) + q(a,b)\otimes\psi(x, y),$$
where
$$\phi: U_{\rho}\times U_{\sigma}\rightarrow U_{\tau}\text{~and~}\psi: U_{\rho}\times U_{\sigma}\rightarrow U_{\tau}$$
is some skew-symmetric and symmetric bilinear maps respectively.
\end{lemma}
Preserving the assumptions and notation of the introduction, we consider a short $SL_2$-structure on a simple Lie algebra $\mathfrak{g}$. In this paragraph, we will derive the basic commutational formulas that arise when commuting two arbitrary elements of isotypic components of the algebra $\mathfrak{g}$.

A commutator on $\mathfrak{g}$ gives $SL_2$-equivariant bilinear maps:
\begin{align*}
\mathfrak{g}_i\times\mathfrak{g}_j\longrightarrow\mathfrak{g}_i\otimes\mathfrak{g}_j&\longrightarrow\mathfrak{g}\qquad(i\neq j),\\
\mathfrak{g}_i\times\mathfrak{g}_i\longrightarrow\textstyle\bigwedge^2\mathfrak{g}_i&\longrightarrow\mathfrak{g}.
\end{align*}
According to the Clebsch-Gordan formula for the tensor product of irreducible representations of the algebra $\mathfrak{sl}_2$ we have:
\begin{equation}\label{gleb}
V_i\otimes V_j\simeq V_{i+j}\oplus V_{i+j-2}\oplus...\oplus V_{|i-j|}.
\end{equation}
Therefore, for commutators of isotypic components of the algebra $\mathfrak{g}$, the following commutational relations are fulfilled:
\begin{equation}\label{kom}
\begin{gathered}
\left[\mathfrak{g}_0, \mathfrak{g}_i\right] \subseteq\mathfrak{g}_i,\qquad i=0,1,2, \\
 \left[\mathfrak{g}_1, \mathfrak{g}_1\right]\subseteq\mathfrak{g}_0\oplus\mathfrak{g}_2, \\
\left[\mathfrak{g}_2, \mathfrak{g}_2\right]\subseteq\mathfrak{g}_0\oplus\mathfrak{g}_2,\\
\left[\mathfrak{g}_1, \mathfrak{g}_2\right]\subseteq \mathfrak{g}_1.
\end{gathered}
\end{equation}
In particular, $\mathfrak{g}_0$ and $\mathfrak{g}_0\oplus\mathfrak{g}_2$ are subalgebras, and a very short $SL_2$-structure is induced on the subalgebra $\mathfrak{g}_0\oplus\mathfrak{g}_2$.

We introduce on the space $\mathbb{C}^2$ $SL_2$-invariant skew-symmetric bilinear form $\langle, \rangle$ according to the following formula: $$\langle u, v\rangle := \det(u, v),\qquad\forall u, v\in \mathbb{C}^2.$$

We also define $SL_2$-invariant symmetric bilinear form on $\mathfrak{sl}_2$ by the rule:
$$(X, Y) := \operatorname{tr}(XY),\qquad\forall X, Y\in \mathfrak{sl}_2.$$

Finally, we define a symmetric bilinear map $S:\mathbb{C}^2\times\mathbb{C}^2\rightarrow\mathfrak{sl}_2$ using the following formula:
$$S(u, v)w = \langle w, u\rangle v + \langle w, v\rangle u.$$
For the completeness of the presentation of the material, we will prove one auxiliary proposition, which is well-known among specialists. We will need it for further reasoning.

Consider an arbitrary simple Lie algebra $\mathfrak{g}$ on which an involutive non-real automorphism $\theta$ is given. The eigenvalues of the automorphism $\theta$ are $1$ and $-1$, so the following decomposition will take place
$$\mathfrak{g} = \mathfrak{h}\oplus\mathfrak{v},$$ where
$$\mathfrak{h} = \{\xi\in\mathfrak{g}: \theta(\xi) = \xi\}$$ and $$\mathfrak{v} = \{\eta\in\mathfrak{g}: \theta(\eta) = -\eta\}.$$
\begin{propos}
The subspace $\mathfrak{h}$ is a Lie subalgebra and its adjoint action on the subspace $\mathfrak{v}$ is faithfull.
\end{propos}
\begin{Proof}
On the basis of the definition of an involutive automorphism, the following relations are fulfilled:
\begin{equation*}
\begin{gathered}
\left[\mathfrak{h}, \mathfrak{h}\right] \subseteq\mathfrak{h},\\
 \left[\mathfrak{h}, \mathfrak{v}\right]\subseteq\mathfrak{v}.
\end{gathered}
\end{equation*}
From where it can be seen that $\mathfrak{h}$ is a Lie subalgebra and the attached one acts on the subspace $\mathfrak{v}$. Let's prove that the given action is faithfull.

To do this, consider the kernel of inefficiency of this action:
$$\mathfrak{k}:= \{\xi\in\mathfrak{h}: [\xi, \eta] = 0,\quad\forall\eta\in\mathfrak{v}\}.$$
Note that since $\theta\neq\operatorname{id}$, then $\mathfrak{h}\neq\mathfrak{g}$. Obviously, $\mathfrak{k}\triangleleft\mathfrak{h}$, and since $\mathfrak{g}= \mathfrak{h}\oplus\mathfrak{v}$ and $\mathfrak{k}$ commutes with $\mathfrak{v}$, then $\mathfrak{k}$ is an ideal of the Lie algebra $\mathfrak{g}$, from which we conclude that $\mathfrak{k} = 0$ or $\mathfrak{k}= \mathfrak{g}$, but by definition $\mathfrak{k}\subset\mathfrak{h}\neq\mathfrak{g}$, so $\mathfrak{k} = 0$. Thus, the adjoint action of $\mathfrak{h}$ on $\mathfrak{v}$ is faithfull.
\end{Proof}

Now let $\mathfrak{g}$ have a short $SL_2$-structure. Note that the action of the element $-\operatorname{id}$ of the group $SL_2$ on the algebra $\mathfrak{g}$ sets an involutive automorphism on it, for which the subalgebra $\mathfrak{g}_0\oplus\mathfrak{g}_2$ is a fixed point space, and $\mathfrak{g}_1$ --- a subspace on which this automorphism acts by multiplying by $-1$. It follows from the sentence proved above that the adjoint action of the algebra $\mathfrak{g}_0\oplus\mathfrak{g}_2$ on the space $\mathfrak{g}_1$ is faithfull. It is follows, that the elements of the space $J_2$, as well as the algebra $\mathfrak{g}_0$, can be identified with their corresponding operators in $J_1$.

From this moment, we will agree to denote the elements of the space $J_1$ with small Latin letters $a, b, c...$, and the elements of the space $J_2$ and the algebra $\mathfrak{g}_0$ we will denote with large Latin letters $A, B, C,...$.

According to lemma \ref{lem1}, as well as the representation of the elements of the space $J_2$ and the algebra $\mathfrak{g}_0$ by operators in $J_1$, it is not difficult to derive the following commutational formulas:
\begin{equation*}
[D, u\otimes a] = u\otimes Da,\qquad\forall D\in \mathfrak{g}_0, u\in \mathbb{C}^2, a\in J_1;
\end{equation*}
\begin{equation*}
[D, X\otimes A] = X\otimes\nu(D, A),\qquad\forall D\in \mathfrak{g}_0, X\in \mathfrak{sl}_2, A\in J_2;
\end{equation*}
\begin{equation*}
[X\otimes A, v\otimes b] = Xv\otimes Ab,\qquad\forall v\in \mathbb{C}^2, X\in\mathfrak{sl}_2, b\in J_1, A\in J_2.
\end{equation*}
Here $\nu: \mathfrak{g}_0\times J_2\rightarrow J_2$ is bilinear map.

To calculate the map $\nu$, it is necessary to consider the action of the operator $\nu(D, A)\in J_2$ on an arbitrary vector $b\in J_1$.
To do this, consider the Jacobi identity for $X\otimes A\in \mathfrak{g}_2,$ $v\otimes b\in\mathfrak{g}_1$ and $D\in\mathfrak{g}_0$:
\begin{equation*}
\bigl[[D, X\otimes A], v\otimes b\bigr]+ \bigl[[X\otimes A, v\otimes b], D\bigr] + \bigl[[v\otimes b , D], X\otimes A\bigr] = 0.
\end{equation*}

Using the commutational relations mentioned above, it is not difficult to conclude that:
$$\nu(D, A) = [D, A],\qquad\forall D\in\mathfrak{g}_0, A\in J_2.$$
Thus, using the formula above, the action of the Lie algebra $\mathfrak{g}_0$ on the space $J_2$ is given.

Commutation of two elements from $\mathfrak{g}_2$ according to the lemma \ref{lem2} is carried out by the formula:
$$
[X\otimes A, Y\otimes B] = [X, Y]\otimes(A\circ B) + \frac{1}{2}(X, Y)\Delta(A, B),\quad\forall X, Y\in \mathfrak{sl}_2, A, B\in J_2.
$$
Here the icon $\circ$ denotes a commutative binary operation that gives the space $J_2$ the structure of a commutative algebra, and $\Delta: J_2\times J_2\rightarrow \mathfrak{g}_0$ is a skew-symmetric bilinear map. The coefficient $\frac{1}{2}$ is taken for the convenience of further calculations.

Thus, we have a vector space $J_1$ on which Lie algebra $\mathfrak{g}_0$ acts, as well as a commutative algebra $J_2$ of linear operators of the space $J_1$ (with respect to multiplication $\circ$). In the following we will show that the space $J_1$ is a symplectic space, the algebra $J_2$ is a subalgebra in the Jordan algebra of all symmetric operators of the space $J_1$, and the Lie algebra $\mathfrak{g}_0$ is a subalgebra in the Lie algebra $\mathfrak{sp}(J_1)$.

Using all the arrangements described above, you can write out a complete list of commutational formulas:
\begin{equation}\label{kom'1}
[D, u\otimes a] = u\otimes Da,\qquad\forall D\in \mathfrak{g}_0, u\in \mathbb{C}^2, a\in J_1;
\end{equation}
\begin{equation}\label{kom'2}
[D, X\otimes A] = X\otimes [D, A],\qquad\forall D\in \mathfrak{g}_0, X\in \mathfrak{sl}_2, A\in J_2;
\end{equation}
\begin{equation}\label{kom'3}
[u\otimes a, v\otimes b] = S(u, v)\otimes\varphi(a, b) + \langle u, v\rangle\delta(a, b),\qquad\forall u, v\in \mathbb{C}^2, a, b\in J_1;
\end{equation}
\begin{equation}\label{kom'4}
[X\otimes A, v\otimes b] = Xv\otimes Ab,\qquad\forall v\in \mathbb{C}^2, X\in\mathfrak{sl}_2, A\in J_2, b\in J_1;
\end{equation}
\begin{equation}\label{kom'5}
[X\otimes A, Y\otimes B] = [X, Y]\otimes(A\circ B) + \dfrac{1}{2}(X, Y)\Delta(A, B),\,\,\forall X, Y\in \mathfrak{sl}_2, A, B\in J_2.
\end{equation}
Here:

$\varphi: J_1\times J_1 \rightarrow J_2$ is a skew-symmetric bilinear map,

$\delta: J_1\times J_1 \rightarrow \mathfrak{g}_0$ is a symmetric bilinear map,

$\Delta: J_2\times J_2 \rightarrow \mathfrak{g}_0$ is a skew-symmetric bilinear map,

$\circ$ is a commutative binary operation on $J_2$.

\subsection{Properties of isotypic components}\label{isprop}

Consider an arbitrary short $SL_2$-structure on a simple Lie algebra $\mathfrak{g}$. Let's write out the Jacobi identities on the Lie algebra $\mathfrak{g}$, from which we derive the identities we need for further reasoning.

For any vector $u\in\mathbb{C}^2$ we will denote by $u^*$ such an element of space $(\mathbb{C}^{2})^*$, that $u^*(v) = \langle v, u\rangle,\quad\forall v\in\mathbb{C}^2$. Then the following formula for map $S$ will be true (here tensors of type (1,1) are interpreted as linear operators):
$$S(u, v) = u\otimes v^* + v\otimes u^*,\qquad\forall u, v\in\mathbb{C}^2.$$
Now we prove one auxiliary lemma simplifying further calculations.

\begin{lemma}\label{lem3}
For any three vectors $u, v, w\in\mathbb{C}^2$ and operator $X\in\mathfrak{sl}_2$ the following equations will be true:
\begin{enumerate}[label=\alph*)]
  \item $\langle Xu, v\rangle = -\langle u, Xv\rangle;$
  \item $u\otimes (Xv)^* = -(u\otimes v^*)\cdot X,\,Xv\otimes u^* = X\cdot(v\otimes u^*);$
  \item $(X, S(u,v)) = 2\langle Xu, v\rangle;$
  \item $\langle u, v\rangle w + \langle v, w\rangle u + \langle w, u\rangle v = 0;$
  \item $u\otimes v^* - v\otimes u^* = \langle u, v\rangle\operatorname{id};$
  \item $[X, S(u, v)] = 2\bigl(S(u, Xv) + \langle u, v\rangle X\bigr);$
  \item $S(Xu, v) - S(u, Xv) = 2\langle u, v\rangle X$.
\end{enumerate}
\end{lemma}
\begin{Proof}

$a)$ It is obvious that $\langle Bu, Bv\rangle = \langle u, v\rangle,\quad\forall B\in SL_2, u, v\in\mathbb{C}^2$. Differentiating this equation in the unit, we get the required one.

$b)$ Since these equations are operator equations, we will check them in relation to an arbitrary vector $w\in\mathbb{C}^2$:

$$(Xv\otimes u^*)w = \langle w, u\rangle Xv = X(\langle w, u\rangle v) = \bigl(X\cdot(v\otimes u^*)\bigr)w;$$

$$(u\otimes (Xv)^*)w = \langle w, Xv\rangle u = -\langle Xw, v\rangle u =\bigl(-(u\otimes v^*)\cdot X\bigr)w.$$
\begin{multline*}
 c)\,(X, S(u,v)) =\operatorname{tr}(X\cdot(u\otimes v^* + v\otimes u^*)) =\operatorname{tr}((u\otimes v^*)\cdot X) + \operatorname{tr}((v\otimes u^*)\cdot X)=\\
= -\operatorname{tr}(u\otimes (Xv)^*)-\operatorname{tr}(v\otimes (Xu)^*) = \langle Xv, u\rangle + \langle Xu, v\rangle = \\
=\langle Xu, v\rangle + \langle Xu, v\rangle = 2\langle Xu, v\rangle.
\end{multline*}

$d)$ The left part of this identity is a trilinear skew-symmetric map on the space $\mathbb{C}^2$. Therefore, it is equal to zero.

$e)$ We prove this operator identity for an arbitrary $w\in\mathbb{C}^2$. We have:
\begin{equation*}
(u\otimes v^* - v\otimes u^*)w =\langle w, v\rangle u - \langle w, u\rangle v = -\langle v, w\rangle u - \langle w, u\rangle v = \langle u, v\rangle w.
\end{equation*}
\begin{multline*}
 f)\,2S(u, Xv)+2\langle u, v\rangle X = 2(u\otimes (Xv)^* +(Xv)\otimes u^*) + 2\langle u, v\rangle X =\\
=-2(u\otimes v^*)\cdot X + 2X\cdot(v\otimes u^*) + 2(u\otimes v^*)\cdot X - 2(v\otimes u^*)\cdot X  = \\
 = [X, 2(v\otimes u^*)] = [X, u\otimes v^* + v\otimes u^* + v\otimes u^* - u\otimes v^*] =\\
 = [X, S(u,v)] + [X, v\otimes u^* - u\otimes v^*]
= [X, S(u,v)] + [X, \langle v, u\rangle\operatorname{id}] = [X, S(u,v)].
\end{multline*}
\begin{multline*}
 g)\,S(Xu, v) - S(u, Xv) = Xu\otimes v^* + v\otimes (Xu)^* - u\otimes (Xv)^* - Xv\otimes u^*=\\=X\cdot(u\otimes v^* - v\otimes u^*) + (u\otimes (v)^* - v\otimes (u)^*)\cdot X=\\=\langle u, v\rangle X  +  \langle u, v\rangle X = 2\langle u, v\rangle X.
\end{multline*}
\end{Proof}

Let us proceed directly to the conclusion of the identities we need for further reasoning. To do this, we will need some facts from the theory of Jordan algebras, a summary of which will be given below. All these facts can be found in\cite{Albert}.
\begin{definition}
An algebra $J$ in which the identities $ab =ba$, $(a^2b)a = a^2(ba)$ are valid for arbitrary $a, b\in J$ is called a Jordan algebra.
\end{definition}
Let $A$ be an arbitrary associative algebra. Vector space $A$ with Jordan multiplication operation:
$$a\circ b = \frac{1}{2}(ab+ba),\qquad\forall a, b\in A,$$
forms the algebra $A^{+}$, which is Jordan. A Jordan algebra embedded in the algebra $A^{+}$ for some associative algebra $A$ is called a {\it special Jordan algebra}.

Denote by $L_c$ the operator of multiplication by the element $c\in J$ of some Jordan algebra $J$. Then transformations of the form $[L_a, L_b]$, where $a, b\in J$ are differentiations of the algebra $J$. Their linear combinations are called {\it internal differentiations}. They form an ideal, which denotes $\mathfrak{inn}(J)$, in the Lie algebra $\mathfrak{der}(J)$ of all differentiations.

\begin{definition}
A Jordan algebra $J$ is called semisimple if a non-degenerate scalar product $(,)$ with the following associativity property is given on the algebra $J$:
$$(ab, c) = (a, bc),\qquad \forall a, b, c\in J.$$
\end{definition}
\begin{definition}
A Jordan algebra $J$ is called simple if it does not contain nontrivial (other than $0$ and $J$) ideals.
\end{definition}
A simple Jordan algebra is semisimple. If the Jordan algebra $J$ is semisimple, then, firstly, it is a direct sum of its simple ideals, secondly, $J$ contains the unit, and, thirdly, $\mathfrak{der}(J) = \mathfrak{inn}(J)$.

The following theorem is also true for semisimple Jordan algebras:
\begin{theorem}\label{semjord}
If $J$ is a finite-dimensional semisimple Jordan algebra over the field $\mathbb{C}$, then the algebra $\mathfrak{der}(J)$ is semisimple.
\end{theorem}
\begin{Proof}
 It follows from \cite[chapter~VIII,~paragraph~4]{Jacobson} that for a semisimple Jordan algebra $J$ over $\mathbb{C}$, the algebra of differentiations $\mathfrak{der}(J)$ is reductive. According to the theory of very short $SL_2$-structures constructed in \cite[section~3]{Vinberg_1}, for a simple Jordan algebra $J$, the algebra of its differentiations $\mathfrak{der}(J)$ coincides with the subspace $\mathfrak{g}_0\subset\mathfrak{g}$, where $\mathfrak{g}$ is a simple Lie algebra with a very short $SL_2$-structure (that is, having the form (\ref{ttk})), which uniquely corresponds to a given Jordan algebra $J$. From the classification of all such structures, it follows that the component $\mathfrak{g}_0$ for a very short $SL_2$-structure is semisimple. For a semisimple Jordan algebra, all differentiations are internal, so the Lie algebra $\mathfrak{der}(J)$ will also be semisimple.
\end{Proof}

Let's return to the consideration of a short $SL_2$-structure on a simple Lie algebra $\mathfrak{g}$. It is clear that a very short $SL_2$-structure is induced on the reductive Lie subalgebra $\mathfrak{g}_0\oplus\mathfrak{g}_2$ of the Lie algebra $\mathfrak{g}$. From the theory of very short $SL_2$-structures, in particular, it follows that the algebra $J_2$ is endowed with the structure of the Jordan algebra and the Lie algebra $\mathfrak{g}_0$ acts on $J_2$ by differentiations.

For convenience, in the future we will call the Jordan algebra $J_2$ \textit{by the Jordan algebra of a short $SL_2$-structure} on $\mathfrak{g}$. Finally, using the Jacobi identity, the multiplication on the algebra $J_2$ is completely determined:
\begin{theorem}\label{th2} Jordan algebra $J_2\subset\mathfrak{gl}(J_1)$ of the short $SL_2$-structure on a simple Lie algebra $\mathfrak{g}$ is special with classical multiplication defined by the formula:
\begin{equation}\label{jordop}
A\circ B = \dfrac{1}{2}(AB + BA),\qquad\forall A, B\in J_2.
\end{equation}
Moreover, the map $\Delta: J_2\times J_2 \rightarrow\mathfrak{g}_0$ satisfies the following relation:
$$\Delta(A, B) = [A, B],\quad \forall A, B\in J_2,$$
where equation above is understood as equation of linear operators on $J_1$.
\end{theorem}
\begin{Proof}
Consider the Jacobi identity for two elements from $\mathfrak{g}_2$ and one element from $\mathfrak{g}_1$  ($w\in\mathbb{C}^2$, $X, Y\in\mathfrak{sl}_2, A, B\in J_2, c\in J_1$):
\begin{equation*}
\bigl[[X\otimes A, Y\otimes B], w\otimes c\bigr] + \bigl[[Y\otimes B, w\otimes c], X\otimes A\bigr] + \bigl[[w\otimes c, X\otimes A], Y\otimes B\bigr] = 0.
\end{equation*}
Using the commutational relations (\ref{kom'1})--(\ref{kom'5}), from this equation we obtain:
\begin{equation}\label{jac1}
YXw\otimes BAc + [X, Y]w\otimes(A\circ B)c + \dfrac{1}{2}(X, Y)w\otimes\Delta(A, B)c - XYw\otimes ABc = 0.
\end{equation}
In equation (\ref{jac1}), we will swap $A$ and $B$, taking advantage of the commutativity of multiplication by $J_2$ and the skew symmetry of $\Delta$. We get:
\begin{equation}\label{jac12}
YXw\otimes ABc + [X, Y]w\otimes(A\circ B)c - \dfrac{1}{2}(X, Y)w\otimes\Delta(A, B)c - XYw\otimes BAc = 0.
\end{equation}
Adding the identity (\ref{jac12}) with (\ref{jac1}), we get:
$$YXw\otimes (BAc + ABc) + 2[X, Y]w\otimes(A\circ B)c - XYw\otimes (BAc + ABc) = 0;$$
$$[Y, X]w\otimes (BAc + ABc) + 2[X, Y]w\otimes(A\circ B)c = 0;$$
$$[Y, X]w\otimes (BAc + ABc - 2(A\circ B)c) = 0,$$
from where we get the first statement of the theorem.

Subtracting from (\ref{jac1}) the identity (\ref{jac12}), we get
$$(X, Y)w\otimes\Delta(A, B)c + YXw\otimes (BAc - ABc) - XYw\otimes (ABc - BAc) = 0;$$
$$(X, Y)w\otimes\Delta(A, B)c - (YXw + XYw)\otimes ([A, B]c)  = 0.$$
It is not difficult to verify that:
$$(XY+YX)w = (X, Y)w,\qquad\forall X,Y\in\mathfrak{sl}_2, w\in \mathbb{C}^2.$$
Hence the second statement of the theorem follows.
\end{Proof}

In connection with the theorem just proved, in the following arguments, wherever the Jordan algebra is used, by default we will assume that this algebra is special with classical multiplication defined by the formula $(\ref{jordop})$.

For simplicity of further reasoning, we will understand by the image of an arbitrary linear map $\psi: U\times V\rightarrow W$ the image of an one-to-one corresponding map from $U\otimes V$ to $W$. It follows from the theory of very short $SL_2$-structures that the image of the $\Delta$ map acts on the Jordan algebra $J_2$ by internal differentiations, and this action is faithfull. Then it follows from the theorem \ref{th2} and the identity (\ref{kom'2}) that this action is calculated by the following formula for arbitrary $A, B, C\in J_2$:
\begin{multline*}
  [[A, B], C] = ABC + CBA  - BAC - CAB=\\=4\bigl(A\circ(B\circ C) - B\circ (A\circ C)\bigr) = 4[L_{A}, L_{B}](C).
\end{multline*}

Thus, we can assume that
\begin{equation*}
\mathfrak{inn}(J_2) = [J_2, J_2] = \bigl\langle[A, B]:\- A,B\in J_2\bigr\rangle.
\end{equation*}

\subsection{Invariant scalar multiplication}\label{invscalmult}

On a simple Lie algebra $\mathfrak{g}$, there is a unique invariant scalar multiplication up to a scalar multiplication. We fix it and in the future we will denote this multiplication with parentheses $(,)$. It follows from Schur's lemma that with respect to this scalar multiplication on $\mathfrak{g}$ vector spaces $\mathfrak{g}_i$, where $i\in\{0, 1, 2\}$, are pairwise orthogonal to each other, that is, the following relations are fulfilled:
\begin{equation}\label{ort}
(\mathfrak{g}_i, \mathfrak{g}_j) = 0,\quad i\neq j.
\end{equation}
Using fixed scalar multiplication on the Lie algebra $\mathfrak{g}$, we introduce scalar multiplication on the space $J_1$ and algebra $J_2$ as follows.
$$(u\otimes a, v\otimes b) = \langle u, v\rangle\alpha(a, b);$$
$$(X\otimes A, Y\otimes B) = \frac{1}{2}(X,Y)\beta(A,B),$$
where $\alpha$ is the skew-symmetric bilinear form on $J_1$ and $\beta$ --- symmetric bilinear form on $J_2$. The coefficient $\frac{1}{2}$ in the last formula is necessary for the convenience of further calculations.

It is obvious that both forms obtained are non-degenerate. Let's call the form $\alpha$ {\it a skew-scalar multiplication on $J_1$}, and the form $\beta$ --- {\it a scalar multiplication on $J_2$}. For brevity, in the future we will denote the form $\beta$ using parentheses $(,)$, and the form $\alpha$ using angle brackets $\langle,\rangle$.

For the introduced scalar and skew-scalar multiplications , the following sentence is true:
\begin{propos}\label{scalth}
Scalar multiplications on the symplectic space $J_1$ and the algebras $\mathfrak{g}_0$ and $J_2$ have the following properties:
\begin{enumerate}[label=\arabic*)]
\item $\bigl(A, \varphi(a, b)\bigr) = \langle Aa, b\rangle = \langle a, Ab\rangle,\qquad\forall a, b\in J_1, A\in J_2$;
\item $\bigl(D, \delta(a, b)\bigr) = \langle Da, b\rangle = -\langle a, Db\rangle,\qquad\forall a, b\in J_1, D\in\mathfrak{g}_0$;
\item $\bigl(D, [A, B]\bigr) = \bigl([D, A], B\bigr) = - \bigl(A, [D, B]\bigr),\qquad\forall A, B\in J_2, D\in\mathfrak{g}_0$;
    \item $(A\circ B, C) = (A, B\circ C),\qquad\forall A, B, C\in J_2$.
\end{enumerate}
\end{propos}
\begin{Proof}

1). Consider $X\in\mathfrak{sl}_2, u, v\in\mathbb{C}^2, a, b\in J_1, A\in J_2$. Then it follows from the invariance of scalar multiplication on the Lie algebra $\mathfrak{g}$ that:
\begin{equation}\label{scal1}
\bigl(X\otimes A, [u\otimes a, v\otimes b]\bigr) = \bigl([X\otimes A, u\otimes a], v\otimes b\bigr).
\end{equation}
Using the relations (\ref{kom'1})--(\ref{kom'5}), we get:
$$\bigl(X\otimes A, [u\otimes a, v\otimes b]\bigr) = \bigl(X\otimes A, S(u, v)\otimes\varphi(a, b)\bigr);$$
$$\bigl([X\otimes A, u\otimes a], v\otimes b\bigr) = \bigl(Xu\otimes Aa, v\otimes b\bigr).$$
From where, using the definition of scalar multiplication introduced above on the spaces $J_1$ and $J_2$, we get that:
$$\dfrac{1}{2}\bigl(X, S(u, v)\bigr)\bigl(A,\varphi(a,b)\bigr) = \langle Xu, v\rangle\langle Aa, b\rangle.$$
According to lemma \ref{lem3} $(X,S(u,v)) = 2\langle Xu, v\rangle$, from which the first equation of the point 1 follows.

The second equation of point 1 easily follows from the first with the replacement of $a$ by $b$ and $b$ by $a$ and the skew-symmetry of the $\varphi$ map.

2). Consider $D\in\mathfrak{g}_0, u, v\in\mathbb{C}^2, a, b\in J_1$. Then it follows from the invariance of scalar multiplication on the Lie algebra $\mathfrak{g}$ that:
\begin{equation}\label{scal2}
\bigl(D, [u\otimes a, v\otimes b]\bigr) = \bigl([D, u\otimes a], v\otimes b\bigr).
\end{equation}
Using the relations (\ref{kom'1})--(\ref{kom'5}) and (\ref{ort}), we get:
$$\bigl(D, [u\otimes a, v\otimes b]\bigr) =\langle u, v\rangle\bigl(D, \delta(a, b)\bigr);$$
$$\bigl([D, u\otimes a], v\otimes b\bigr) = \bigl(u\otimes Da, v\otimes b\bigr).$$
From where, using the definition of skew-scalar multiplication on the space $J_1$ introduced above, we get
$$\langle u, v\rangle\bigl(D, \delta(a, b)\bigr) = \langle u, v\rangle\langle Da, b\rangle.$$
Dividing the resulting equation by $\langle u, v\rangle$, we get the first equation of point 2. The second equation of point 2 easily follows from the first with the replacement of $a$ by $b$ and $b$ by $a$ and the symmetry of $\delta$.

3). Consider $D\in\mathfrak{g}_0, X, Y\in\mathfrak{sl}_2, A, B\in J_2$. Then it follows from the invariance of scalar multiplication on the Lie algebra $\mathfrak{g}$ that:
\begin{equation}\label{scal3}
\bigl(D, [X\otimes A, Y\otimes B]\bigr) = \bigl([D, X\otimes A], Y\otimes B\bigr).
\end{equation}
Using the relations (\ref{kom'1})--(\ref{kom'5}) and the theorem \ref{th2}, from (\ref{scal3}) we get:
$$\dfrac{1}{2}\bigl(X, Y\bigr)\bigl(D, [A, B]\bigr) = \bigl(X\otimes[D, A], Y\otimes B\bigr) = \dfrac{1}{2}\bigl(X, Y\bigr)\bigl([D, A], B\bigr).$$

Hence we get the first equation of point 3. The second equation of point 3 easily follows from the first with the replacement of $A$ by $B$ and $B$ by $A$.

4). Consider $X, Y, Z\in\mathfrak{sl}_2, A, B, C\in J_2$. Then it follows from the invariance of scalar multiplication on the Lie algebra $\mathfrak{g}$ that:
\begin{equation}\label{scal4}
\bigl(X\otimes A, [Y\otimes B, Z\otimes C]\bigr) = \bigl([X\otimes A, Y\otimes B], Z\otimes C\bigr).
\end{equation}
Using the relations (\ref{kom'1})--(\ref{kom'5}), from (\ref{scal4}) we get:
$$\dfrac{1}{2}\bigl(X, [Y, Z]\bigr)\bigl(A, B\circ C\bigr) = \dfrac{1}{2}\bigl([X, Y], Z\bigr)\bigl(A\circ B, C\bigr)$$
Since $\bigl(X, [Y,Z]\bigr) = \bigl([X, Y], Z\bigr)$, then from the equation above you can get the equation of point 4.
\end{Proof}

From the proposition \ref{scalth} and the theorem \ref{th2} in particular, it follows that $J_2$ is a semisimple Jordan algebra and $\mathfrak{der}(J_2) = \mathfrak{inn}(J_2)=[J_2, J_2]\subset\mathfrak{g}_0$. It also follows from the proposition \ref{scalth} that the algebra $J_2$ consists of symmetric operators of the symplectic space $J_1$, while $\mathfrak{g}_0\subset \mathfrak{sp}(J_1)$.

\subsection{Short $SL_2$-structures and $\mathbb{Z}$-graduations}\label{zgrad}
Consider an arbitrary short $SL_2$-structure on a simple Lie algebra $\mathfrak{g}$. Using the map $\operatorname{d}\Phi$ Lie algebra $\mathfrak{sl}_2$ is embedded in $\mathfrak{der}(\mathfrak{g}) = \mathfrak{inn}(\mathfrak{g})\simeq\mathfrak{g}$. Denote by $e, f, h$ the basic elements of the algebra $\operatorname{d}\Phi(\mathfrak{sl}_2)$, which satisfy the following relations:
\begin{equation}\label{sl2}
[e, f] = h,\quad [h, e]= 2e,\quad [h, f] = -2f.
\end{equation}
Let
$$\tilde{e}= \operatorname{d}\Phi^{-1}(e),\quad\tilde{f}= \operatorname{d}\Phi^{-1}(f),\quad\tilde{h}= \operatorname{d}\Phi^{-1}(h).$$
Consider the operator $\operatorname{ad}(h)$. Since the dimensions of the irreducible representations of $\Phi$ do not exceed 3, then the eigenvalues of the operator $\operatorname{ad}(h)$ do not exceed modulo 2. Then, with respect to the operator $\operatorname{ad}(h)$, the algebra $\mathfrak{g}$ decomposes into a direct sum of its own subspaces:
\begin{equation}\label{dec}
\mathfrak{g} = \mathfrak{g}^{-2}\oplus\mathfrak{g}^{-1}\oplus\mathfrak{g}^{0}\oplus\mathfrak{g}^{1}\oplus\mathfrak{g}^2,\quad\mathfrak{g}^{k} = \{\xi\in\mathfrak{g}: [h,\xi] = k\xi\}.
\end{equation}
It is obvious, that
$$\mathfrak{g}^0 = \mathfrak{g}_0\oplus(\tilde{h}\otimes J_2),\quad\mathfrak{g}^{-1}= e_{-1}\otimes J_1,\quad\mathfrak{g}^1 =e_{1}\otimes J_1,\quad\mathfrak{g}^{-2} = \tilde{f}\otimes J_2,\quad\mathfrak{g}^{2}=\tilde{e}\otimes J_2,$$
where $\{e_{-1}, e_{1}\}$ is the proper basis of the space $\mathbb{C}^2$ with respect to the operator $\tilde{h}$.

The decomposition (\ref{dec}) is $\mathbb{Z}$-grading, that is, $[\mathfrak{g}^{k},\mathfrak{g}^{l}]\subset\mathfrak{g}^{k+l}$, where $\mathfrak{g}^m = 0$ for $|m|> 2$.

Without limiting of generality, we can assume that the element $h$ belongs to a fixed Cartan subalgebra $\mathfrak{t}$ of the algebra $\mathfrak{g}$ and even a positive Weyl chamber with respect to some choice of positive roots.

Thus, if $\Pi = \{\alpha_1,...,\alpha_{n}\}$ is the system of simple roots of algebra $\mathfrak{g}$, then $$\alpha_{i}(h)\geqslant0,\qquad\forall i=\overline{1,n}.$$
Denote $p_i = \alpha_i(h)\in\mathbb{Z_{+}}.$ It is clear that $\mathbb{Z}$-grading (\ref{dec}) is completely determined by the set $\{p_1,...,p_n\}$. Each of the subspaces $\mathfrak{g}^k$ is the sum of the root subspaces $\mathfrak{g}_{\alpha}$ corresponding to the roots $\alpha = k_1\alpha_1+ ... + k_n\alpha_n$ with the condition:
\begin{equation}\label{rooteq1}
k_1p_1 + \dots + k_np_n = k
\end{equation}
and Cartan subalgebras $\mathfrak t$ in the case of $k=0$.

Denote by $\alpha_{\operatorname{\operatorname{max}}}$ the highest root of the Lie algebra $\mathfrak{g}$. Since the maximum eigenvalue of the operator $\operatorname{ad}(h)$ is 2, then $\alpha_{\operatorname{\operatorname{max}}}(h) = 2$. If $\alpha_{\operatorname{\operatorname{max}}} = l_1\alpha_1 + \dots + l_n\alpha_n$, then from the condition that $\alpha_{\operatorname{\operatorname{max}}}(h) = 2$, we get that
\begin{equation*}
l_1p_1 + \dots + l_np_n = 2.
\end{equation*}

Consider the subspace $\mathfrak{g}^0$. Obviously, $\mathfrak{g}^{0}$ is a Lie subalgebra.
The following decomposition takes place: $$\mathfrak{g}^{0} = \tilde{\mathfrak{g}}^{0}\oplus\langle h\rangle,$$ where $\tilde{\mathfrak{g}}^{0}$ is orthogonal complement to $h$. Consider the representation $\tilde{\mathfrak{g}}^{0}$ on the space $\mathfrak{g}^{2}$, which is a constraint of the adjoint representation.
Denote $$\Pi_0 = \{\alpha_i\in\Pi: p_i = 0\}.$$
It is clear that $\Pi_0$ is a system of simple roots of the semisimple part of $\mathfrak{\tilde{g}}^{0}$.

Consider the subspace $\mathfrak{g}^k$ for $k\neq 0$. Decompose it into a direct sum of subspaces $\mathfrak{g}^k_{(\nu)}$, where each $\mathfrak{g}^k_{(\nu)}$ is the sum of the root subspaces $\mathfrak{g}_{\alpha}$ corresponding to roots with fixed coefficients for simple roots $\alpha_i\not\in\Pi_0$ (and with the condition (\ref{rooteq1})).

Obviously, each of the subspaces $\mathfrak{g}^k_{(\nu)}$ is invariant with respect to $\mathfrak{g}^0$. The following statement is true, the proof of which can be found in \cite[chapter~3,~$\S$~3,~paragraph~3.5]{Vinberg_3}:
\begin{theorem}
Representation of $\mathfrak{g}^0$ in each $\mathfrak{g}^k_{(\nu)}$ is irreducible.
\end{theorem}
It immediately follows from this theorem that the representation $\tilde{\mathfrak{g}}^0$ in each $\mathfrak{g}^k_{(\nu)}$ is irreducible.

Note that in our case the representation $\tilde{\mathfrak{g}}^0:\mathfrak{g}^2$ is always irreducible, since, generally speaking, only three cases are possible:
\begin{enumerate}
\item there is exactly one simple root $\alpha_i$ with a non-zero value $p_i=2$, and in the decomposition of the highest root by simple roots, this root stands with a coefficient of ~1;
\item there are exactly two simple roots $\alpha_i$ and $\alpha_j$ with nonzero values $p_i = p_j = 1$, in the decomposition of the highest root by simple roots, the coefficients for these roots are equal to 1;
\item there is exactly one simple root $\alpha_i$ with a non-zero value $p_i=1$, and in the decomposition of the highest root by simple roots, this simple root stands with a coefficient of 2.
\end{enumerate}
Case 1 is not implemented, because with it $\mathfrak{g}^1 = \mathfrak{g}^{-1}=0$, which is not the case for short $SL_2$-structures. In case 3, there is one simple root, so $\mathfrak{g}^2 = \mathfrak{g}^2_{(\nu)}$, which means the representation is irreducible. In case 2, in the decomposition of the highest root by simple roots, the coefficients for roots with nonzero values by $h$ are equal to $1$, which means that in all positive roots, the coefficients for these roots do not exceed $1$, that is, $\mathfrak{g}^2 = \mathfrak{g}^2_{(\nu)}$, and, again, the representation is irreducible.
\subsection{Simplicity of the Jordan algebra of a short $SL_2$-structure}\label{simjord}
\begin{theorem}\label{simp}
Let a short $SL_2$-structure be given on a simple Lie algebra $\mathfrak{g}$ and a $J_2$ is the Jordan algebra of this structure. Then $J_2$ is simple.
\end{theorem}
\begin{Proof}
This fact is easy to prove for very short $SL_2$-structures. Indeed, in this case, the Lie algebra $\mathfrak{g}$ has the decomposition $\mathfrak{g} = \mathfrak{g}_0\oplus\mathfrak{sl}_2\otimes J_2$, where $J_2$ is a Jordan algebra, and $\mathfrak{g}_0 = \mathfrak{der}(J_2)=\mathfrak{inn}(J_2)$. If the Jordan algebra $J_2$ has any nontrivial ideal $I$, then the subspace $\mathfrak{g}_0\oplus\mathfrak{sl}_2\otimes I$ is a nontrivial ideal in the algebra $\mathfrak{g}$, which follows easily from the commutation formulas.

In the case of short $SL_2$-structures, such a direct method of proof is inapplicable, since when commuting with itself the summand $\mathbb{C}^2\otimes J_1$, one summand arises from $\mathfrak{sl}_2\otimes J_2$, which prevents the above-described method of constructing an ideal in $\mathfrak{g}$. However, the case of short $SL_2$-structures can be reduced to the already proved case of very short $SL_2$-structures.

To do this, note that if we restrict a short $SL_2$-structure on a simple Lie algebra $\mathfrak{g}$ to a Lie subalgebra $\mathfrak{h} = \mathfrak{g}_0\oplus\mathfrak{g}_2$, we get a very short $SL_2$-structure given by on the Lie algebra $\mathfrak{h}$. This Lie algebra, generally speaking, is not simple, but it is reductive, due to the non-degeneracy of scalar multiplication. Since the semisimple element $h$ (like the whole algebra $\mathfrak{sl}_2$) trivially acts on the center of the algebra $\mathfrak{h}$, then replacing the algebra with its commutator $[\mathfrak{h}, \mathfrak{h}]$, we can assume that the algebra $\mathfrak{h}$ is semisimple.

For convenience, we will choose a system of simple roots for $\mathfrak{h}$ as follows. Obviously, the root system of the algebra $\mathfrak{h}$ is the roots of the algebra $\mathfrak{g}$, taking the values $0$, $2$ or $-2$ on the element $h$. We will consider positive roots of the algebra $\mathfrak{g}$ taking the value 0 on the element $h$, as well as the roots of the algebra $\mathfrak{g}$ taking the value $-2$ on $h$ as the positive roots of the algebra $\mathfrak{h}$. Then it is clear that the simple roots of the semisimple algebra $\mathfrak{h}$ are all the simple roots of the Lie algebra $\mathfrak{g}$ taking zero value on the element $h$, as well as the lowest root of the algebra $\mathfrak{g}$. It is obvious that, generally speaking, the Lie algebra $\mathfrak{h}$ contains several simple components, but exactly one of them contains the subspace $\mathfrak{sl}_2\otimes J_2$ (namely, the one whose simple root is the lowest root of the algebra $\mathfrak{g}$). Therefore, we can consider the algebra $\mathfrak{h}$ simple, which completely reduces the question to the case of very short $SL_2$-structures, and for it it has already been proved that the Jordan algebra $J_2$ is simple.
\end{Proof}

Consider the Lie algebra $\mathfrak{g}_0$. In connection with the commutation relations (\ref{kom}), the subspace $[J_2, J_2]\subset\mathfrak{g}_0$ is the ideal of the algebra $\mathfrak{g}_0$. Denote by $\mathfrak{i}_0$ the orthogonal complement to $[J_2, J_2]$ in $\mathfrak{g}_0$ relative to the scalar multiplication introduced above. The subspace $\mathfrak{i}_0\subset\mathfrak{g}_0$ is also an ideal of the algebra $\mathfrak{g}_0$ and, due to the non-degeneracy of scalar multiplication by $\mathfrak{g}_0$, the sum of the dimensions of the ideals $\mathfrak{i}_0$ and $[J_2, J_2]$ is equal to the dimension of the Lie algebra $\mathfrak{g}_0$.

Consider the action of the ideal $\mathfrak{i}_0\triangleleft\mathfrak{g}_0$ on the Jordan algebra $J_2$. Using the proposition \ref{scalth} proved above, we have $\forall A, B\in J_2, D\in \mathfrak{i}_0$:
\begin{equation}{\label{eqsc1}}
\bigl([D, A], B\bigr) = \bigl(D, [A, B]\bigr) = 0.
\end{equation}
Therefore, $[D, A] = 0,\-\forall A\in J_2$, that is, $D$ is an element of the kernel of the action $\mathfrak{g}_0$ on $J_2$. Conversely, if $[D, A] = 0,\-\forall A\in J_2$, then it follows from the equation (\ref{eqsc1}) that $\bigl(D, [A,B]\bigr) = 0,\-\forall A, B\in J_2$, that is, $D\in\mathfrak{i}_0$. From which it immediately follows that $\mathfrak{i}_0$ is exactly the kernel of the action of $\mathfrak{g}_0$ on $J_2$. From the fact that the Lie algebra $\mathfrak{g}_0$ acts on $J_2$ by commutation, and multiplication on $J_2$ is classical, it is not difficult to deduce that the Lie algebra $\mathfrak{g}_0$ acts on the Jordan algebra $J_2$ by differentiations. Then, since $\mathfrak{i}_0$ is the kernel of this action, then $\mathfrak{g}_0/\mathfrak{i}_0\subset\mathfrak{der}(J_2)$. But from the simplicity of the algebra $J_2$ it follows that
$$\mathfrak{der}(J_2) = [J_2, J_2],$$
that is
$$\mathfrak{g}_0/\mathfrak{i}_0\subset[J_2, J_2].$$
Then there is a decomposition of the Lie algebra $\mathfrak{g}_0$ into a direct sum of two ideals:
\begin{equation}\label{g0dec}
\mathfrak{g}_0 = \mathfrak{i}_0\oplus [J_2, J_2].
\end{equation}

From the classification of all short $SL_2$-structures described in the section \ref{classif} of this article, it will follow, in particular, that, generally speaking, $\mathfrak{i}_0\neq 0$, that is, the Lie algebra $\mathfrak{g}_0$ is not completely defined Jordan algebra $J_2$, which distinguishes the case of short $SL_2$-structures from very short ones. Moreover, the Lie algebra $\mathfrak{g}_0$ is not completely defined even by the pair $(J_1; J_2)$, which will also be seen from the classification of all short $SL_2$-structures.

Note that from the simplicity of the Jordan Lie algebra $J_2$ it follows that $J_2$ contains the unit. Moreover, the following theorem is true:
\begin{theorem}
Let a short $SL_2$-structure be given on a simple Lie algebra $\mathfrak{g}$, $J_1$ is a symplectic space, and $J_2$ is a Jordan algebra of this structure. Denote the unit element $J_2$ by $\mathbb{I}$. Then $\mathbb{I} = \operatorname{id},$ that is $$\mathbb{I}a=a,\quad\forall a\in J_1.$$
\end{theorem}
\begin{Proof}

To begin with, we will deduce the identities we need to prove. Consider the Jacobi identity for arbitrary $X\in\mathfrak{sl}_2, u,v\in\mathbb{C}^2, A\in J_2, a, b\in J_1$:
\begin{equation*}
\bigl[[X\otimes A, u\otimes a], v\otimes b\bigr] + \bigl[[u\otimes a, v\otimes b], X\otimes A\bigr] + \bigl[[v\otimes b, X\otimes A], u\otimes a\bigr] = 0.
\end{equation*}
Using the equations (\ref{kom'1})--(\ref{kom'5}), we get:
\begin{multline}\label{keyjac}
S(Xu, v)\otimes\varphi(Aa, b) + \langle Xu, v\rangle\delta(Aa, b) + \\ + \bigl[S(u, v), X\bigr]\otimes\bigl(\varphi(a, b)\circ A\bigr) + \frac{1}{2}\bigl(S(u, v), X\bigr)\bigl[\varphi(a, b), A\bigr] + \\ + \langle u, v\rangle X\otimes\bigl[\delta(a, b), A\bigr] - S(Xv, u)\otimes\varphi(Ab, a) - \langle Xv, u\rangle\delta(Ab, a) = 0.
\end{multline}
Where from
\begin{equation*}
\langle Xu, v\rangle\delta(Aa, b) + \frac{1}{2}\bigl(S(u, v), X\bigr)\bigl[\varphi(a, b), A\bigr] - \langle Xv, u\rangle\delta(Ab, a) = 0.
\end{equation*}
Applying the lemma \ref{lem3} to the obtained equation, we obtain the first of the identities we need:
\begin{equation}\label{iid1}
\bigl[A, \varphi(a, b)\bigr] = \delta(Aa, b) - \delta(a, Ab).
\end{equation}
Further, it also follows from (\ref{keyjac}) that:
\begin{multline}\label{id2'}
S(Xu, v)\otimes\varphi(Aa, b) + \bigl[S(u, v), X\bigr]\otimes\bigl(\varphi(a, b)\circ A\bigr) + \\ + \langle u, v\rangle X\otimes\bigl[\delta(a, b), A\bigr] - S(Xv, u)\otimes\varphi(Ab, a) = 0.
\end{multline}
In the resulting equation, we will swap the elements $a$ and $b$ and use the skew symmetry of $\varphi$ and the symmetry of $\delta$. We get:
\begin{multline}\label{id2"}
S(Xu, v)\otimes\varphi(Ab, a) - \bigl[S(u, v), X\bigr]\otimes\bigl(\varphi(a,b)\circ A\bigr) + \\ + \langle u, v\rangle X\otimes\bigl[\delta(a, b), A\bigr] - S(Xv, u)\otimes\varphi(Aa, b) = 0.
\end{multline}
Adding the equations (\ref{id2'}) and (\ref{id2"}), we get:
\begin{equation*}
\bigl(S(Xu, v) - S(Xv, u)\bigr)\otimes \bigl(\varphi(Ab, a) + \varphi(Aa, b)\bigr)+ 2\langle u, v\rangle X\otimes\bigl[\delta(a, b), A\bigr] = 0.
\end{equation*}
From the lemma \ref{lem3} it follows that $S(Xu,v) - S(Xv, u) = 2\langle u, v\rangle X$, from where we get the second of the identities we need:
\begin{equation}\label{id2}
\bigl[A, \delta(a, b)\bigr] = \varphi(Aa, b) - \varphi(a, Ab).
\end{equation}
Consider the Jacobi identity for arbitrary $D\in\mathfrak{g}_0, u, v\in\mathbb{C}^2, a, b\in J_1$:
\begin{equation*}
\bigl[[D, u\otimes a], v\otimes b\bigr] + \bigl[[u\otimes a, v\otimes b], D\bigr] + \bigl[[v\otimes b, D], u\otimes a\bigr] = 0.
\end{equation*}
Using the equations (\ref{kom'1})--(\ref{kom'5}), we get:
\begin{multline*}
S(u, v)\otimes\varphi(Da, u) + \langle u, v\rangle\delta(Da, b) -\\ - S(u, v)\otimes \bigl[D, \varphi(a, b)\bigr] - \langle u, v\rangle\bigl[D, \delta(a, b)\bigr] + \\ + S(u, v)\otimes\varphi(a, Db) + \langle u, v\rangle\delta(a, Db) = 0.
\end{multline*}
From the equation above we have:
\begin{equation*}
\langle u, v\rangle\delta(Da, b) - \langle u, v\rangle\bigl[D, \delta(a, b)\bigr] + \langle u, v\rangle\delta(a, Db) = 0.
\end{equation*}
Dividing the resulting equation by $\langle u, v\rangle$, we get:
\begin{equation}\label{id3}
\bigl[D, \delta(a, b)\bigr] = \delta(Da, b) + \delta(a, Db).
\end{equation}
Now we have everything we need to prove the statement of the theorem.

Consider the operator $\mathbb{I}\in J_2$. Since $\mathbb{I}^2=\mathbb{I}$, then $J_1 = J_1^0\oplus J_1^1$, where
$$J_1^0 = \operatorname{Ker}\mathbb{I},\qquad J_1^1 = \operatorname{Im}\mathbb{I}.$$
We will prove that $J_1^0 = 0$, from which the statement of the theorem will follow.

It is clear that $J_2 = J_2 \circ\mathbb{I}$, which means that the entire Jordan algebra $J_2$ trivially acts on $J_1^0$.

We prove that $J_1^0$ is an invariant subspace for the entire Lie algebra $\mathfrak{g}_0$. Due to the triviality of the action of the Jordan algebra $J_2$ on $J_1^0$, the commutator $[J_2, J_2]\subset\mathfrak{g}_0$ also trivially acts on $J_1^0$, therefore, by virtue of the decomposition (\ref{g0dec}), it is sufficient to prove that $J_1^0$ is invariant with respect to the action of the ideal $\mathfrak{i}_0$. We have for $\forall D\in\mathfrak{i}_0, a\in J_1^0$
$$\mathbb{I}(Da)=[\mathbb{I},D]a+D(\mathbb{I}a)=0 \implies Da\in J_1^0.$$
This proves the invariance of $J_1^0$ under the action of $\mathfrak{g}_0$.
Denote
$$\mathfrak{g}_0^1 = \bigl\langle\delta(a, b):\, a\in J_1^0, b\in J_1\bigr\rangle.$$

Consider $\mathfrak{i} = \mathfrak{g}_0^1\oplus\bigl(\mathbb{C}^2\otimes J_1^0\bigr)$. We will prove that $\mathfrak{i}\triangleleft\mathfrak{g}$, whence, by virtue of simplicity $\mathfrak{g}$, it will follow that $J_1^0 = 0$.
Note that
$$\bigl(A, \varphi(a, b)\bigr) = \langle Aa, b\rangle = 0,\qquad\forall A\in J_2, a\in J_1^0, b\in J_1.$$
Therefore,
\begin{equation}\label{eqphi}
\varphi(a, b) = 0,\quad\forall a\in J_1^0, b\in J_1.
\end{equation}
From the invariance of the subspace $J_1^0$ with respect to the action of the Lie algebra $\mathfrak{g}_0$, using the equation (\ref{id3}), it is not difficult to obtain that the subspace $\mathfrak{i}$ is invariant with respect to commutation with elements from $\mathfrak{g}_0$. Also, using the obtained equations (\ref{eqphi}) and (\ref{id2}) and the triviality of the action $J_2$ on $J_1^0$, it can be obtained that the subspace $\mathfrak{i}$ commutes with $\mathfrak{g}_2$. It remains to prove that the subspace $\mathfrak{i}$ is invariant with respect to commutation with elements from $\mathfrak{g}_1$. Since the equation (\ref{eqphi}) is fulfilled, it is enough to show that the operators from $\mathfrak{g}_0^1$ map an arbitrary element from $J_1$ to $J_1^0$.

Applying the equation (\ref{iid1}) to $a\in J_1^0$ and $A = \mathbb{I}$ and considering (\ref{eqphi}), we can get that
$$\delta(a, \mathbb{I}b) = \delta(\mathbb{I}a, b) = 0,\quad\forall a\in J_1^0, b\in J_1.$$
Thus we get that
\begin{equation}\label{eqdel}
\delta(a, b) = 0,\quad\forall a\in J_1^0, b\in J_1^1.
\end{equation}
Therefore, it is sufficient to consider the operators $\delta(a,b)$ for $a,b\in J_1^0$.

We show that $\delta(a,b)c=0$, $\forall a,b\in J_1^0$, $c\in J_1^1$. Then it follows from the invariance of $J_1^0$ with respect to the action of $\mathfrak{g}_0$ that the operators from $\mathfrak{g}_0^1$ map an arbitrary element from $J_1$ to $J_1^0$.

Consider the Jacobi identity for arbitrary $u,v,w\in\mathbb{C}^2, a, b\in J_1^0, c\in J_1^1$:
\begin{equation*}
\bigl[[u\otimes a, v\otimes b], w\otimes c\bigr] + \bigl[[v\otimes b, w\otimes c], u\otimes a\bigr] + \bigl[[w\otimes c, u\otimes a], v\otimes b\bigr] = 0.
\end{equation*}
Using the equations (\ref{kom'1})--(\ref{kom'5}) as well as the equations (\ref{eqphi}) and (\ref{eqdel}), we get:
$$\bigl[[u\otimes a, v\otimes b], w\otimes c\bigr] = \langle u, v\rangle w\otimes \delta(a, b)c,$$
$$\bigl[[v\otimes b, w\otimes c], u\otimes a\bigr] = \bigl[[w\otimes c, u\otimes a], v\otimes b\bigr] = 0.$$
Thus, it follows from the Jacobi identity above that $\delta(a,b)c = 0,\,\forall a, b\in J_1^0, c\in J_1^1$.

Hence, $\mathfrak{i}\triangleleft\mathfrak{g}$, and thus the statement of the theorem is proved.
\end{Proof}

Note that from the proof of the theorem above it follows that for an arbitrary short $SL_2$-structure, the map $\Delta: J_1\times J_1\rightarrow\mathfrak{r}_0$, like the map $\varphi: J_1\times J_1\rightarrow J_2$, is equivariant with respect to the action $\mathfrak{r}_0$ (identity (\ref{id3}), and the identities (\ref{iid1}) and (\ref{id2}) are also executed. Since decomposition (\ref{g0dec}) takes place, then
$$\delta = \delta_0 + \delta_c,\quad\delta_0 = \pi_{0,0}\delta,\quad\delta_c = \pi_{0,c}\delta,$$
where $\pi_{0,0}: \mathfrak{g}_0\rightarrow\mathfrak{i}_0$ and $\pi_{0,c}: \mathfrak{g}_0\rightarrow [J_2, J_2]$ are orthogonal projections on $\mathfrak{i}_0$ and $[J_2, J_2]$ respectively.

Since the map $\delta$ is $\mathfrak{g}_0$-equivariant, it is also true for $\delta_0$. Since the identity (\ref{iid1}) is fulfilled, the left part of which belongs to $[J_2, J_2]$, then
$$\delta_0(A a, b) = \delta_0(a, A b),\quad\forall A\in J_2, a, b\in J_1.$$

\section{Maximal short $SL_2$-structure}\label{maxsl2str}

The purpose of this section is to build an example of a short $SL_2$-structure that will be maximal in the sense that for this short $SL_2$-structure, the components $J_2$ and $\mathfrak{g}_0$ are the largest in inclusion for a given component $J_1$.

\subsection{Construction of the maximal short $SL_2$-structure}\label{maxsl2build}

In this subsection, we will consider the algebra $\mathfrak{g} = \mathfrak{so}_{4n+1}$ in a basis in which its elements are matrices of size $(4n+1)\times(4n+1)$, skew-symmetric with respect to the side diagonal. As an invariant scalar product on $\mathfrak{so}_{4n+1}$, we fix the following:
$$(A, B) = \dfrac{1}{2}\operatorname{tr}(AB),\quad\forall A, B\in\mathfrak{so}_{4n+1}.$$

Consider on $\mathfrak{so}_{4n+1}$ a short $SL_2$-structure given by embedding the Lie algebra $\mathfrak{sl}_2$ in $\mathfrak{so}_{4n+1}$, in which the basis elements $e,f$ and $h$ of algebra $\mathfrak{sl}_2\subset\mathfrak{so}_{4n+1}$ satisfying the relations (\ref{sl2})
are matrices of the following form:
$$e = \begin{pmatrix}
0&I_{2n}\\
0&0
\end{pmatrix},\quad f = \begin{pmatrix}
0&0\\
I_{2n}&0
\end{pmatrix},\quad h = \operatorname{diag}\{\underbrace{1,..,1}_{2n}, 0, \underbrace{-1,..,-1}_{2n}\}.$$
Here $I_{2 n}$ denotes the matrix of size $2n\times 2n$ of the following form:
$$I_{2n} =  \operatorname{diag}\{\underbrace{1,..,1}_{n}, \underbrace{-1,..,-1}_{n}\}.$$
Note that $$\mathfrak{so}_{4n+1}=\mathfrak{so}_{4n}\oplus\mathbb{C}^{4n}.$$
This decomposition is easily obtained by considering an arbitrary element of the algebra $\mathfrak{so}_{4n+1}$ in the above basis, which has the following form:
$$\begin{pmatrix}
&&\overline{a}_1&&\\
&$\Huge\text{$C$}$&\overline{a}_2&$\Huge\text{$A$}$&\\
\overline{b}_{2}^{\mathsf{s}}&-\overline{b}_{1}^{\mathsf{s}}&0&-\overline{a}_{2}^{\mathsf{s}}&-\overline{a}_{1}^{\mathsf{s}}\\
&$\Huge\text{$B$}$&\overline{b}_{1}&$\Huge\text{$-C^{\mathsf{s}}$}$&\\
&&-\overline{b}_{2}&&
\end{pmatrix},$$
where $$\overline{a}_1 = \begin{pmatrix}
a_{1}\\
\vdots\\
a_{n}
\end{pmatrix},\quad\overline{a}_2 = \begin{pmatrix}
    a_{n+1}\\
\vdots\\
a_{2n}
\end{pmatrix},\quad\overline{b}_1 = \begin{pmatrix}
b_{1}\\
\vdots\\
b_{n}
\end{pmatrix},\quad\overline{b}_2 = \begin{pmatrix}
b_{n+1}\\
\vdots\\
b_{2n}
\end{pmatrix}.$$
Here $A, B$ are matrices of size $2 n\times2n$ skew-symmetric with respect to the side diagonal, and $C$ is a matrix of size $2n\times 2 n$, and $C^{\mathsf{s}}$ is a matrix transposed to $C$ relative to the side diagonal. A direct summand of the form $\mathfrak{so}_{4 n}$ forms matrices in which $a_i = b_i = 0,\,\forall i=\overline{1, 2 n}$. A direct summand of the form $\mathbb{C}^{4 n}$ forms matrices in which the blocks $A, B, C$ and $C^{\mathbf{s}}$ are zero, and the matrices from the first direct summand naturally act on the second direct summand.

The following equalities also take place: $$\mathfrak{so}_{4n}=\mathfrak{g}_0\oplus\mathfrak{g}_2,\quad\mathbb{C}^{4n}=\mathfrak{g}_1.$$ Next, the quadratic space $\mathbb{C}^{4 n}$ decomposes into the tensor product of the symplectic spaces $\mathbb{C}^2$ and $\mathbb{C}^{2n}=J_1$. Namely, vector from $\mathbb{C}^{4 n}$, which is a matrix of the type "cross" (that is, the matrix $(4 n+1)\times(4 n+1)$, whose nonzero elements belong only to the central column and the central row):
$$\begin{pmatrix}
&&\overline{a}_1&&\\
&&\overline{a}_2&&\\
\overline{b}_{2}^{\mathsf{s}}&-\overline{b}_{1}^{\mathsf{s}}&0&-\overline{a}_{2}^{\mathsf{s}}&-\overline{a}_{1}^{\mathsf{s}}\\
&&\overline{b}_{1}&&\\
&&-\overline{b}_{2}&&
\end{pmatrix},$$
is mapped to an element of the space $\mathbb{C}^2\otimes J_1$, having the form:
$$e_1\otimes(a_1,...,a_{2n}) + e_{-1}\otimes(b_1,...,b_{2n}).$$
Considering two arbitrary elements of the space $\mathfrak{g}_1$ and multiplying them scalar, it is not difficult to deduce that the matrix of skew-scalar multiplication on the space $J_1$ has the form:
$$\Omega = \begin{pmatrix}
&\Omega_n\\
-\Omega_n&
\end{pmatrix}.$$
Here, $\Omega_n$ denotes a matrix of size $n\times n$ of the following form:
$$
\Omega_n =\begin{pmatrix}
&&1\\
&\rotatebox{80}{$\ddots$}&\\
1&&
\end{pmatrix}.
$$
The space of skew-symmetric operators in $\mathbb{C}^{4n}$ decomposes into the direct sum of two tensor products: the space of skew-symmetric operators in $\mathbb{C}^2$ and the space of symmetric operators in $J_1$, and, conversely, the space of symmetric operators in $\mathbb{C}^2$ and space of skew-symmetric operators in $J_1$. The first of these products is $\mathfrak{sl}_2\otimes J_2=\mathfrak{g}_2$, and since the space of symmetric operators in $\mathbb{C}^2$ is obviously one-dimensional, we can assume that the second term is the algebra $\mathfrak{sp}(J_1)=\mathfrak{g}_0$. Let's denote the Jordan algebra of symmetric operators of the symplectic space $J_1$ by $\mathfrak{sym}(J_1)$. So, $J_2 = \mathfrak{sum}(J_1)$.

Based on the form of the matrix of skew-scalar multiplication by $J_1$, it is not difficult to obtain a matrix description of the algebra $\mathfrak{sum}(J_1)$ in the basis of the space $J_1$, in which the skew-scalar multiplication has the matrix $\Omega$ specified above. Namely:
$$ \mathfrak{sym}(J_1) = J_2 = \left\{\begin{pmatrix}
X&Y\\
Z&X^{\mathsf{s}}
\end{pmatrix}\in\mathfrak{gl}_{2n}: X, Y, Z\in\mathfrak{gl}_n,\,Y^{\mathsf{s}} = -Y,\,Z^{\mathsf{s}} = -Z\right\}.$$
Similarly, the matrix description of the algebra $\mathfrak{sp}(J_1)$ in the same basis has the following form:
$$ \mathfrak{sp}(J_1) = \mathfrak{g}_0 = \left\{\begin{pmatrix}
X&Y\\
Z&-X^{\mathsf{s}}
\end{pmatrix}\in\mathfrak{gl}_{2n}: X, Y, Z\in\mathfrak{gl}_n,\,Y^{\mathsf{s}} = Y,\,Z^{\mathsf{s}} = Z\right\}.$$
We show explicitly how tensor products of elements of the Jordan algebra $J_2=\mathfrak{sym}(J_1)$ on the elements of the Lie algebra $\mathfrak{sl}_2$ correspond to the elements of the Lie algebra $\mathfrak{so}_{4n+1}$. Let's reinterpret through $e, f$ and $h$ the standard basis of the abstract Lie algebra $\mathfrak{sl}_2$ satisfying the relations (\ref{sl2}).

The matrix $M\in\mathfrak{sym}(J_1)$, having the form:
$$M = \begin{pmatrix}
X&Y\\
Z&X^{\mathsf{s}}
\end{pmatrix},\text{~where~} X, Y, Z\in\mathfrak{gl}_n,\,Y^{\mathsf{s}} = -Y,\,Z^{\mathsf{s}} = -Z,$$
and the basis element $h\in\mathfrak{sl}_2$ correspond to the matrix $\hat{M}_{h}\in h\times J_2\subset\mathfrak{so}_{4 n+1}$, which has the form:
$$\hat{M}_{h} = \begin{pmatrix}
&&0&&\\
&$\Huge \text{$M$}$&\vdots&&\\
0&\dots&0&\dots&0\\
&&\vdots&$\Huge\text{$-M^{\mathsf{s}}$}$&\\
&&0&&
\end{pmatrix}.$$
Here and further, there are zeros in the empty spaces in the matrices.

Similarly, the matrix $M$ and the basis element $e\in\mathfrak{sl}_2$ correspond to the matrix $\hat{M}_{e}\in e\times J_2\subset\mathfrak{so}_{4n+1}$, which has the form:
$$\hat{M}_{e} = \begin{pmatrix}
&&0&&\\
&&\vdots&$\Huge\text{$M_{c}$}$&\\
0&\dots&0&\dots&0\\
&&\vdots&&\\
&&0&&
\end{pmatrix},$$
where
$$M_{c} = M \cdot I_{2n} = \begin{pmatrix}
X&-Y\\
Z&-X^{\mathsf{s}}
\end{pmatrix}.$$

Finally, the matrix $M$ and the basis element $f\in\mathfrak{sl}_2$ correspond to the matrix $\hat{M}_{f}\in f\otimes J_2\subset\mathfrak{so}_{4n+1}$, which has the form:
$$\hat{M}_{f} = \begin{pmatrix}
&&0&&\\
&&\vdots&&\\
0&\dots&0&\dots&0\\
&$\Huge\text{$M_{r}$}$&\vdots&&\\
&&0&&
\end{pmatrix},$$
where
$$M_{r}= I_{2n} \cdot M = \begin{pmatrix}
X&Y\\
-Z&-X^{\mathsf{s}}
\end{pmatrix}.$$
Explicit embedding $\mathfrak{g}_0=\mathfrak{sp}(J_1)\subset\mathfrak{so}_{4n+1}$ is arranged as follows. The matrix $N\in\mathfrak{sp}(J_1)$, having the form:
$$N = \begin{pmatrix}
X&Y\\
Z&-Z^{\mathsf{s}}
\end{pmatrix},\text{~where~} X, Y, Z\in\mathfrak{gl}_n,\,Y^{\mathsf{s}} = Y,\,Z^{\mathsf{s}} = Z,$$
corresponds to the matrix $\hat{N}_{0}\in\mathfrak{so}_{4n+1}$, which has the form:
$$\hat{N}_{0} =\begin{pmatrix}
&&0&&\\
&$\Huge\text{$N$}$&\vdots&&\\
0&\dots&0&\dots&0\\
&&\vdots&$\Huge\text{$-N^{\mathsf{s}}$}$&\\
&&0&&
\end{pmatrix}.$$

Let's denote the scalar multiplication by $J_2$ using $(,)_2$. Considering two arbitrary elements of the space $\mathfrak{g}_2$ and multiplying them scalar, it is not difficult to obtain that scalar multiplication on the algebra $J_2\simeq\mathfrak{sym}(J_1)$ has the form:
$$(A,B)_2 = \operatorname{tr}(AB),\qquad\forall A, B\in\mathfrak{sym}(J_1). $$
Denote scalar multiplication by $\mathfrak{g}_0$ with $(,)_0$. Acting similarly with the elements of the algebra $\mathfrak{g}_0\simeq\mathfrak{sp}(J_1)$, we obtain that the scalar multiplication on it has the form:
$$(D_1, D_2)_0 = \operatorname{tr}(D_1D_2),\qquad\forall D_1, D_2\in \mathfrak{sp}(J_1). $$
Thus, summing up the above, we have that for the considered short $SL_2$-structure it is fulfilled that:
$$J_1= \mathbb{C}^{2n},\quad J_2 = \mathfrak{sym}(J_1),\quad\mathfrak{g}_0=\mathfrak{sp}(J_1),$$
moreover the scalar multiplication on $\mathfrak{g}_0$ and on $J_2$ are  standard scalar multiplications of linear operators on $J_2$.

The constructed short $SL_2$-structure will be called {\it maximal short $SL_2$-structure}. Since for any short $SL_2$-structure with symplectic space $J_1$, inclusions of $J_2\subset\mathfrak{sym}(J_1)$ are performed and $\mathfrak{g}_0\subset\mathfrak{sp}(J_1)$, then for the maximal short $SL_2$-structures with symplectic space $J_1$, the algebras $J_2$ and $\mathfrak{g}_0$ are the largest in inclusion, which justifies such a name.

\subsection{Analysis of the maximal short $SL_2$-structure}\label{maxsl2anal}

Consider on the symplectic space $J_1$ linear operators $R(a, b)$ of rank 1 $(a, b\in J_1)$ defined by the formula:
$$R(a,b)c = \langle c, a\rangle b,\qquad\forall a, b, c\in J_1.$$
Obviously, the introduced operator $R$ has the following property: $$R^*(a, b) = -R(b, a),\qquad\forall a, b\in J_1,$$ where $R^*(a,b): J_1\rightarrow J_1$ --- operator conjugate to the operator $R(a, b)$.
\newline
\newline
Now we will enter the following maps on $J_1$:

1. Map $\varphi_m: J_1\times J_1\rightarrow\mathfrak{gl}(J_1)$ by the formula:
$$\varphi_m(a, b) = \dfrac{1}{2}\bigl(R(b, a) - R(a, b)\bigr).$$

2. Map $\delta_m: J_1\times J_1\rightarrow\mathfrak{gl}(J_1)$ by the formula:
$$\delta_m(a, b) = \dfrac{1}{2}\bigl(R(b, a) +R(a,b)\bigr).$$
It is clear that for arbitrary $a, b\in J_1$ operator $\delta_m(a,b)$ is skew-symmetric, and the operator $\varphi_m(a,b)$ is symmetric with respect to the skew-symmetric scalar multiplication of the space $J_1$, that is:
$$\varphi_m(a, b)\in\mathfrak{sym}(J_1),\quad \delta_m(a, b)\in\mathfrak{sp}(J_1),\qquad\forall a, b\in J_1.$$
From what was proved in the previous paragraph for the maximal short $SL_2$-structure $J_2 = \mathfrak{sym}(J_1)$ and $\mathfrak{g}_0 = \mathfrak{sp}(J_1)$. From the definition of the maps $\delta_m$ and $\varphi_m$, it is easy to conclude that for arbitrary $a, b\in J_1, A\in J_2, D\in\mathfrak{g}_0$:
$$\bigl(\varphi_{m}(a, b), A\bigr)_2 = \langle Aa, b\rangle;$$
$$\bigl(\delta_{m}(a, b), D\bigr)_0 = \langle Da, b\rangle.$$

Therefore, by virtue of the proposition \ref{scalth}, the non-degeneracy of scalar multiplication on $J_2$ and $\mathfrak{g}_0$ and the above equations, we obtain that the maps $\delta_m$ and $\varphi_m$ coincide with the maps $\delta$ and $\varphi$ arising in commutational formulas for the maximal short $SL_2$-structure.

Consider an arbitrary short $SL_2$-structure on the algebra $\mathfrak{g}$, for which $\mathfrak{g}_1 = \mathbb{C}^2\otimes J_1$. Let $\mathfrak{g}_0$ be the isotypic component of the trivial representation of $\mathfrak{sl}_2$, and $J_2$ be the Jordan algebra of this $SL_2$-structure. Due to the simplicity of the Jordan algebra $J_2$ restriction of the invariant scalar multiplication $(,)_2$ of the Jordan algebra $\mathfrak{sym}(J_1)$ on a simple subalgebra $J_2$ is non-degenerate, which means it is proportional to the invariant scalar multiplication on $J_2$ obtained from scalar multiplication on $\mathfrak{g}$. Thus, scalar multiplication on $J_2$ is proportional to the trace of the product of operators on $J_1$. Choosing a suitable coefficient for scalar multiplication on $\mathfrak{g}$, we can assume that the scalar multiplication on $J_2$ has the form
\begin{equation}\label{scalequat2}
(A, B) = (A, B)_2 = \operatorname{tr}(AB),\quad\forall A, B\in J_2.
\end{equation}
Let's use the equality 3 of proposition \ref{scalth}. For arbitrary $D\in\mathfrak{g}_0$ and $A, B\in J_2$ we have:
$$(D, [A,B]) = ([D, A], B) = \operatorname{tr}([D,A]B) = \operatorname{tr}(D[A,B]).$$
Thus, on the subspace $[J_2, J_2]\subset\mathfrak{g}_0$, scalar multiplication is defined by the formula
\begin{equation}\label{scalequat0}
(D_1, D_2) = (D_1, D_2)_0 = \operatorname{tr}(D_1D_2),\quad\forall D_1, D_2\in [J_2, J_2].
\end{equation}
Note that due to the decomposition (\ref{g0dec}) and the non-degeneracy of scalar multiplication on the Lie algebra $\mathfrak{g}_0$, the scalar multiplication calculated above is non-degenerate on $[J_2, J_2]$. For convenience, we will assume everywhere that on a simple Jordan algebra $J_2$, as well as on the algebra of its differentiations $\mathfrak{der}(J_2) = \mathfrak{inn}(J_2 )= [J_2, J_2]$ scalar multiplication is set by default using the formulas (\ref{scalequat2}) and (\ref{scalequat0}), respectively.

\begin{theorem}
For any short $SL_2$-structure with symplectic space $J_1$, the maps $\varphi: J_1\times J_1\rightarrow J_2$ and $\delta: J_1\times J_1\rightarrow\mathfrak{g}_0$ arising in commutation formulas are surjective (as linear maps from $J_1\otimes J_1$), with $\varphi = \pi_{2}\varphi_m,\quad \delta_c = \pi_{c}\delta_m$, where $\pi_{2}: \mathfrak{sym}(J_1)\rightarrow J_2$ and $\pi_{c}:\mathfrak{sp}(J_1)\rightarrow [J_2, J_2]$ are orthogonal projections on the algebras $J_2$ and $[J_2, J_2]$ respectively. For the maximal short $SL_2$-structure $\varphi=\varphi_m$ and $\delta= \delta_m$.
\end{theorem}
\begin{Proof}

For the maximal short $SL_2$-structure of equations $\varphi=\varphi_m$ and $\delta= \delta_m$ have already been proven.

Consider an arbitrary short $SL_2$-structure and the subspace $\mathfrak{g}_1\oplus[\mathfrak{g}_1, \mathfrak{g}_1]$. From the relations (\ref{kom}) it follows that $\mathfrak{g}_1\oplus[\mathfrak{g}_1, \mathfrak{g}_1]\lhd\mathfrak{g}$, from where $\mathfrak{g}_1\oplus[\mathfrak{g}_1, \mathfrak{g}_1] = \mathfrak{g}$. Thus, in connection with the commutational formula (\ref{kom'3}), the following relations take place:
\begin{equation*}
\mathfrak{g}_0 = \left<\delta(a,b):\,a, b\in J_1\right>,\qquad
\mathfrak{g}_2 = \left<\varphi(a,b):\,a, b\in J_1\right>.
\end{equation*}
Hence, the maps $\delta$ and $\varphi$ are surjective.
Acting similarly to the case of the maximal short $SL_2$-structure, for an arbitrary short $SL_2$-structure, according to the proposition $\ref{scalth}$, we obtain that $\forall a, b\in J_1, A, B\in J_2$:
$$(A, \varphi(a, b))_2 = \langle Aa, b\rangle = (A, \varphi_m(a,b))_2,$$
$$([A, B], \delta_c(a, b))_0 = ([A, B], \delta(a, b))_0 = \langle [A, B]a, b\rangle = ([A,B], \delta_m(a, b))_0.$$
Due to the non-degeneracy of scalar multiplications on $J_2$ and $[J_2, J_2]$, the map $\varphi$ coincides with the map $\pi_2\varphi_{m}$, and the map $\delta_c$ is the same as the map $\pi_{c}\delta_{m}$.
\end{Proof}

\section{Short $SL_2$-structures and symplectic Lie-Jordan structures}\label{simpliejord}

\subsection{Short $SL_2$-structures and curvature of symmetric space}\label{levichev}

Let's step back for a moment from considering short $SL_2$-structures.

Let $\mathfrak{h}$ be a reductive Lie algebra and $R:\mathfrak {h}\rightarrow\mathfrak{gl}(\mathfrak{v})$ be faithfull orthogonal linear representation on linear space $\mathfrak{v}$. Let's define the commutators of elements from $\mathfrak{h}$ with elements from $\mathfrak{v}$ using the representation of $R$ by the formula
$$[\xi, x] = R(\xi)x = -[x, \xi],\qquad\forall\xi\in\mathfrak{h},\quad\forall x\in\mathfrak{v},$$
and the commutators of two elements of $\mathfrak{v}$ as elements from $\mathfrak{h}$
satisfying the condition
\begin{equation}\label{komscal}
\bigl(\xi,[x,y]\bigr)=\bigl([\xi,x],y\bigr),\qquad\forall\xi\in\mathfrak{h},\quad\forall x, y\in\mathfrak{v},
\end{equation}
where parentheses denote $\mathfrak{h}$-invariant scalar multiplications in $\mathfrak{h}$ and $\mathfrak{v}$.

On the space $\mathfrak{v}$ we define the 4-linear form
$$K(x, y, z, u) = \bigl([x, y], [z, u]\bigr),\qquad\forall x, y, z, u\in\mathfrak{v}.$$
The following result goes back to E. Cartan and B. Kostant (see, for example, \cite[section~1]{Conlon} and \cite[section~1]{Panyushev}).
\begin{theorem}\label{tenth}
The operation defined above in the space $\mathfrak{h}\oplus\mathfrak{v}$ sets the structure of a $\mathbb{Z}_2$-graded Lie algebra on it if and only if the form $K(x, y, z, u)$ has the symmetry of the curvature tensor, that, in turn, is equivalent to fulfilling the Bianchi identity:
\begin{equation}\label{bian}
K(x, y, z, u) + K(y, z, x, u) + K(z, x, y, u) = 0.
\end{equation}
\end{theorem}
To prove it , it is enough to check the truth of the Jacobi identity in the following four cases:
\begin{enumerate}
\item for elements $\xi, \eta, \zeta\in\mathfrak{h}$;
\item for elements $\xi, \eta\in\mathfrak{h}, x\in\mathfrak{v}$;
\item for elements $\xi\in\mathfrak{h}, x, y\in\mathfrak{v}$;
\item for elements $x, y, z\in\mathfrak{v}$.
\end{enumerate}

In case 1, the Jacobi identity holds because $\mathfrak{h}$ is a Lie algebra. In case 2, the Jacobi identity follows from the fact that $R$ is a representation of the Lie algebra $\mathfrak{h}$. Case 3 can be reduced to case 2 by scalar multiplying the left side of the proved Jacobi identity by an arbitrary element $\kappa\in\mathfrak{h}$ and using scalar multiplication invariance. In the latter case, the Jacobi identity is equivalent to the identity (\ref{bian}).

The Bianchi identity in the theorem above does not appear randomly. The fact is that if $\mathfrak{g} = \mathfrak{h}\oplus\mathfrak{v}$ is $\mathbb{Z}_2$-graded Lie algebra, then the space $G/H$, where $G$ is such a complex group that $\operatorname{Lie} G = \mathfrak{g}$, and $H$ is a Lie subgroup such that $\operatorname{Lie}H= \mathfrak{h}$, is a complex symmetric space on which there is a $G$-invariant metric tensor, induced by scalar multiplication on $\mathfrak{v}$, and symmetric Levi-Civita connectivity preserving it. The Riemann tensor of the symmetric space $G/H$ at the base point $eH$ has the form:
$$\rho(x,y) = -\operatorname{ad}([x,y])|_{\mathfrak{v}},\qquad\forall x, y\in\mathfrak{v},$$
where $\mathfrak{v}$ is identified with tangent space $T_{eH}(G/H)$. Therefore, the 4-line form $K$ is, up to the sign, a covariant curvature tensor on the space $G/H$ at the point $eH$. A more detailed analysis of this construction can be found in \cite[section~1]{Conlon}.

Let's return to the consideration of short $SL_2$-structures. Let's put $\mathfrak{h} = \mathfrak{g}_0\oplus\mathfrak{g}_2, \mathfrak{v}=\mathfrak{g}_1$, and as a representation of $R$ consider the adjoint representation. As previously noted, the adjoint representation of $\mathfrak{h}$ on $\mathfrak{v}$ is faithfull, and due to the invariance of scalar multiplication on $\mathfrak{g}$ is also orthogonal.

Also, due to the properties of invariant scalar multiplication on the Lie algebra $\mathfrak{g}$, the commutation operation satisfies the condition (\ref{komscal}).

For the specified $\mathfrak{h}$ and $\mathfrak{v}$, consider the 4-line form $K$. By the \ref{tenth} theorem, the form $K$ will satisfy the identity (\ref{bian}). In connection with the equalities (\ref{isot}), the 4-linear form $K$ can be rewritten in terms of the space $J_1$. Namely, let's put $x = u_1\otimes a, y = u_2\otimes b, z = u_3\otimes c, u = u_4\otimes d$, where $u_i\in\mathbb{C}^2,\,\forall i=\overline{1, 4},\,\forall a, b, c, d\in J_1$. Then, using the equations (\ref{kom'1})--(\ref{kom'5}), we get:
\begin{equation*}
K(x, y, z, u) = \bigl(S(u_1, u_2)\otimes \varphi(a, b) + \langle u_1, u_2\rangle\delta(a, b), S(u_3, u_4)\otimes \varphi(c, d) + \langle u_3, u_4\rangle\delta(c, d)\bigr).
\end{equation*}
From where, using the definition of scalar multiplication and the lemma \ref{lem3}, we get:
\begin{multline}
    K(x, y, z, u) = \bigl(\langle u_3, u_1\rangle\langle u_2, u_4\rangle + \langle u_3, u_2\rangle\langle u_1, u_4\rangle\bigr)\bigl(\varphi(a, b), \varphi(c, d)\bigr) + \\
    + \langle u_1, u_2\rangle\langle u_3, u_4\rangle\bigl(\delta(a, b), \delta(c, d)\bigr).
\end{multline}
Now the equation (\ref{bian}) can be rewritten as:
 \begin{multline}\label{bian2}
  \langle u_1, u_2\rangle\langle u_3, u_4\rangle\Bigl(\bigl(\varphi(b, c), \varphi(a, d)\bigr) + \bigl(\varphi(a, c), \varphi(b, d)\bigr) + \bigl(\delta(a, b), \delta(c, d)\bigl)\Bigr) +\\
  + \langle u_2, u_3\rangle\langle u_1, u_4\rangle\Bigl(\bigl(\varphi(b, a), \varphi(c, d)\bigr) + \bigl(\varphi(c, a), \varphi(b, d)\bigr) + \bigl(\delta(b, c), \delta(a, d)\bigr)\Bigr) +\\
  + \langle u_3, u_1\rangle\langle u_2, u_4\rangle\Bigl(\bigl(\varphi(a, b), \varphi(c, d)\bigr) + \bigl(\varphi(c, b), \varphi(a, d)\bigr) + \bigl(\delta(c, a), \delta(b, d)\bigr)\Bigr) = 0.
 \end{multline}
Let's denote
 $$\overline{F}(a, b, c, d) = \bigl(\delta(a, b), \delta(c, d)\bigr) + \bigl(\varphi(b, c), \varphi(a, d)\bigr) + \bigl(\varphi(a, c), \varphi(b, d)\bigr).$$
Then the equation (\ref{bian2}) can be rewritten in the form:
 \begin{multline}\label{jac3'}
 \langle u_1, u_2\rangle\langle u_3, u_4\rangle \overline{F}(a, b, c, d) + \langle u_2, u_3\rangle\langle u_1, u_4\rangle \overline{F}(b, c, a, d) + \\ + \langle u_3, u_1\rangle\langle u_2, u_4\rangle \overline{F}(c, a, b, d) = 0.
 \end{multline}

We prove that the identity (\ref{jac3'}) is equivalent to the symmetry of the form $\overline{F}$. To do this, use the well-known Plucker identity for $u_1, u_2, u_3, u_4\in\mathbb{C}^2$:
$$\langle u_1, u_2\rangle\langle u_3, u_4\rangle +\langle u_2, u_3\rangle\langle u_1, u_4\rangle + \langle u_3, u_1\rangle\langle u_2, u_4\rangle = 0.$$

First, (\ref{jac3'}) follows from the symmetry of $\overline{F}$ by virtue of the Plucker identity. Secondly, the definition of $\overline{F}$ immediately implies the symmetry of $\overline{F}$ with respect to the permutation of the first two arguments, as well as with respect to the permutation of the last two arguments. Therefore, to prove the inverse implication, it is enough to make sure that $\overline{F}$ is symmetric with respect to the cyclic permutation of the first three arguments. To do this, substitute in (\ref{jac3'}) $u_1=u_3=e_1, u_2=u_4=e_{-1}$. Then $\langle u_1, u_2\rangle\langle u_3, u_4\rangle = - \langle u_2, u_3\rangle\langle u_1, u_4\rangle = 1$, and $\langle u_3, u_1\rangle\langle u_2, u_4\rangle = 0$, from where it will follow symmetry of $\overline{F}$ with respect to cyclic permutations of the first three arguments.

As an example, we calculate the form $\overline{F}_m$ for the maximum short $SL_2$\-structure. We have:
\begin{multline*}
  \overline{F}_m(a, b, c, d) = \bigl(\delta_m(a, b), \delta_m(c, d)\bigr) + \bigl(\varphi_m(b, c), \varphi_m(a, d)\bigr) + \bigl(\varphi_m(a, c), \varphi_m(b, d)\bigr) =\\
  = \bigl\langle\delta_m(a, b)c + \varphi_m(b, c)a + \varphi_m(a, c)b , d\bigr\rangle =\\
  =\frac{1}{2} \Bigl\langle-\bigl(\langle a, c\rangle b + \langle b, c\rangle a\bigr) + \bigl(\langle b, a\rangle c - \langle c, a\rangle b\bigr) + \bigl(\langle a, b\rangle c - \langle c, b\rangle a\bigr) ,d\Bigr\rangle = 0.
\end{multline*}

In general, we will use the proposition \ref{scalth}:
\begin{multline*}\overline{F}(a, b, c, d) = \bigl(\delta(a, b), \delta(c, d)\bigr) + \bigl(\varphi(b, c), \varphi(a, d)\bigr) + \bigl(\varphi(a, c), \varphi(b, d)\bigr) =\\= \bigl(\delta(a, b)c + \varphi(b, c)a + \varphi(a, c)b, d\bigr).\end{multline*}
Let's denote $$F(a, b, c) = \delta(a, b)c + \varphi(b, c)a + \varphi(a, c)b.$$
Thus, using the formula above, the trilinear map $$F: J_1\times J_1\times J_1\rightarrow J_1$$ is defined.
Since for an arbitrary short $SL_2$-structure, as shown above, the form $\overline{F}$ is symmetric, then the map $F$ is also symmetric.

\subsection{The main theorem}\label{maintheor}
Now it's time to formulate and prove the main theorem of this work. For the convenience of the formulation, we introduce the following definitions.

Let $J_1$ be a symplectic space, $J_2\subset\mathfrak{sym}(J_1)$ be a semisimple Jordan subalgebra whose unit is the identity operator of the space $J_1$ and $\mathfrak{g}_0\subset\mathfrak{sp}(J_1)$ is a reductive Lie subalgebra with $[J_2, J_2]\subset\mathfrak{g}_0$ and $[\mathfrak{g}_0, J_2]\subset J_2$. Then the following decomposition takes place
$$\mathfrak{g}_0 = [J_2, J_2]\oplus \mathfrak{i}_0,$$
where $\mathfrak{i}_0$ is the kernel of the adjoint action $\mathfrak{g}_0$ on $J_2$. Since the Jordan algebra $J_2$ is semisimple, then by the theorem \ref{semjord} the Lie algebra $[J_2, J_2]$ is also semisimple. We assume that a non-degenerate invariant scalar multiplication equal to the trace of the product of operators is given on algebras $\mathfrak{sp}(J_1)$ and $\mathfrak{sym}(J_1)$, as well as on their subalgebras $[J_2, J_2]$ and $J_2$, respectively.

Let be a bilinear symmetric $\mathfrak{g}_0$-equivariant surjective map $\delta_0: J_1\times J_1\rightarrow\mathfrak{i}_0$ satisfying the following identity
\begin{equation}\label{delta0}
\delta_0(Aa, b) = \delta_0(a, Ab),\quad\forall A\in J_2, a, b\in J_1.
\end{equation}
Let's set maps $\delta: J_1\times J_1\rightarrow\mathfrak{g}_0$ and $\varphi: J_1\times J_1\rightarrow J_2$ using formulas
\begin{equation}\label{deltaphi}
\delta = \delta_0 + \delta_c,\quad\delta_c = \pi_{c}\delta_m,\quad\varphi = \pi_2\varphi_m,
\end{equation}
where $\pi_{c}: \mathfrak{sp}(J_1)\rightarrow [J_2, J_2]$ --- orthogonal projection on $[J_2, J_2]$, and $\pi_2: \mathfrak{sym}(J_1)\rightarrow J_2$ --- orthogonal projection on $J_2$.
\begin{definition}
The quad $\left(J_1; J_2; \mathfrak{g}_0; \delta_0\right)$ is called a symplectic Lie-Jordan structure if the trilinear map $F: J_1\times J_1\times J_1\rightarrow J_1$, given by the formula
$$F(a, b, c) = \delta(a, b)c + \varphi(b, c)a + \varphi(a, c)b$$
is symmetric.
\end{definition}
\begin{remark}
It is easy to notice that the maps $\delta_c$ (and hence $\delta$) and $\varphi$  are  $\mathfrak{g}_0$-equivariant. Indeed, $\delta_c = \pi_{c}\delta_m,\,\varphi = \pi_2\varphi_m$, with the maps $\delta_m$ and $\varphi_m$ are $\mathfrak{sp}(J_1)$-equivariant, and the projections $\pi_{c}$ and $\pi_2$  are $\mathfrak{g}_0$-equivariant because they are  orthogonal projections on $\mathfrak{g}_0$-invariant subspaces with respect to invariant scalar multiplication. The map $\delta_c$ is surjective as a composition of the surjective maps $\pi_c$ and $\delta_m$. With the adjoint representation $\mathfrak{g}_0$ on itself, the irreducible terms included in $\mathfrak{i}_0$ and $[J_2, J_2]$ are pairwise non-isomorphic to each other. Therefore, since the maps $\delta_0$ and $\delta_c$ are $\mathfrak{g}_0$-equivariant, then the map $\delta = \delta_0+\delta_c$ will also be surjective as the sum of two surjective homomorphisms of $\mathfrak{g}_0$-modules.
\end{remark}
\begin{definition}
The symplectic Lie-Jordan structure $\left(J_1; J_2; \mathfrak{g}_0; \delta_0\right)$ is called simple if the Jordan algebra $J_2$ is simple.
\end{definition}
\begin{theorem}\label{finth}
There is a one-to-one correspondence between simple Lie-Jordan symplectic structures and simple Lie algebras with a short $SL_2$-structure. Namely, a simple Lie-Jordan symplectic structure $\left(J_1; J_2; \mathfrak{g}_0; \delta_0\right)$ corresponds to the simple Lie algebra $\mathfrak{g}$ with a short $SL_2$-structure, which have the form
\begin{equation}\label{declie}
\mathfrak{g} = \mathfrak{g}_0\oplus\left(\mathbb{C}^2\otimes J_1\right)\oplus\left(\mathfrak{sl}_2\otimes J_2\right),
\end{equation}
where the commutation operation on $ \mathfrak{g}$ is defined by the following formulas ($D\in \mathfrak{g}_0, u,v\in \mathbb{C}^2, a, b\in J_1, X, Y\in\mathfrak{sl}_2, A, B\in J_2$):
\begin{equation}\label{kom''1}
[D, u\otimes a] = u\otimes Da,
\end{equation}
\begin{equation}\label{kom''2}
[D, X\otimes A] = X\otimes [D, A],
\end{equation}
\begin{equation}\label{kom''3}
[u\otimes a, v\otimes b] = S(u, v)\otimes\varphi(a, b) + \langle u, v\rangle\delta(a, b),
\end{equation}
\begin{equation}\label{kom''4}
[ X\otimes B, u\otimes a] = Xu\otimes Ba,
\end{equation}
\begin{equation}\label{kom''5}
[X\otimes A, Y\otimes B] = [X, Y]\otimes(A\circ B) + \dfrac{1}{2}(X, Y)[A, B],
\end{equation}
Here the maps $\delta: J_1\times J_1\rightarrow\mathfrak{g}_0$ and $\varphi: J_1\times J_1\rightarrow J_2$ are given using formulas (\ref{deltaphi}), and $\circ$ denotes the multiplication operation in the Jordan algebra $J_2$ (in formulas (\ref{kom''1}), (\ref{kom''2}) and (\ref{kom''4}), when commuting elements in a different order, a minus sign is implied before the expression on the right).
\end{theorem}
\begin{Proof}

To begin with, we note that for each simple Lie algebra with a short $SL_2$-structure, as shown above, it is possible to construct a simple symplectic Lie-Jordan structure. It remains to prove that a simple Lie-Jordan symplectic structure $(J_1; J_2; \mathfrak{g}_0; \delta_0)$ can be used to construct a simple Lie algebra $\mathfrak{g}$ with a short $SL_2$-structure. Denote $\mathfrak{g}_1 = \mathbb{C}^2\otimes J_1,\,\mathfrak{g}_2= \mathfrak{sl}_2\otimes J_2$. Due to the fact that the formulas (\ref{kom''2}) and (\ref{kom''5}) are identical to the corresponding formulas for a very short $SL_2$-structure, the subspace $\mathfrak{h} = \mathfrak{g}_0\oplus\mathfrak{g}_2$ is endowed with the structure of the Lie algebra.

Consider the space $\mathfrak{g} = \mathfrak{h}\oplus\mathfrak{g}_1$. Let's define a commutation operation on it using the formulas (\ref{kom''1})--(\ref{kom''5}). We prove that with the help of these formulas, the space $\mathfrak{g}$ is endowed with a Lie algebra structure. To do this, it is sufficient to prove that the Jacobi identity is true in the following cases:
\begin{enumerate}
\item for elements $\xi, \eta, \zeta\in\mathfrak{h}$;
\item for elements $\xi, \eta\in\mathfrak{h}, x\in\mathfrak{g}_1$;
\item for elements $\xi\in\mathfrak{h}, x, y\in\mathfrak{g}_1$;
\item for elements $x, y, z\in\mathfrak{g}_1$.
\end{enumerate}
In case 1, the Jacobi identity is fulfilled due to the fact that $\mathfrak{h}$ is a Lie algebra. Using the formulas (\ref{kom''1}) and (\ref{kom''4}), we can determine the linear map $R: \mathfrak{h}\rightarrow\mathfrak{gl}(\mathfrak{g}_1)$ of the Lie algebra $\mathfrak{h}$ in the space of linear operators on $\mathfrak{g}_1$. This map will be a linear representation of the Lie algebra $\mathfrak{h}$. Indeed, we prove the defining relation of the linear representation of the Lie algebra for $R$, namely
\begin{equation}\label{repr}
  R([\xi, \eta])x = [R(\xi), R(\eta)]x = R(\xi)R(\eta)x - R(\eta)R(\xi)x,\quad\forall \xi,\eta\in\mathfrak{h},\,x\in\mathfrak{g}_1.
\end{equation}
Let's say $x = u\otimes a\in\mathfrak{g}_1$. Due to linearity, it is sufficient to prove equality (\ref{repr}) for the following three cases:
\begin{enumerate}
  \item $\xi,\eta\in\mathfrak{g}_0$
  \item $\xi\in\mathfrak{g}_0,\eta\in\mathfrak{g}_2$
  \item $\xi,\eta\in\mathfrak{g}_2$
\end{enumerate}
In the case 1 let's denote $\xi = D_1, \eta = D_2\in\mathfrak{g}_0$. We have:
\begin{equation*}
  R([\xi, \eta])x = u\otimes [D_1, D_2]a = u\otimes D_1D_2a - u\otimes D_2D_1a = [R(\xi), R(\eta)]x.
\end{equation*}
 In the case 2 let's denote $\xi = D\in\mathfrak{g}_0, \eta = X\otimes A\in\mathfrak{g}_2$. We have:
\begin{equation*}
  R([\xi, \eta])x = Xu\otimes [D, A]a = [D, Xu\otimes Aa] - [X\otimes A, u\otimes Da] = [R(\xi), R(\eta)]x.
\end{equation*}
In the case 3 let's denote $\xi = X\otimes A, \eta = Y\otimes B\in\mathfrak{g}_2$. We have:
\begin{multline*}
  R([\xi, \eta])x = [X, Y]u\otimes (A\circ B)a + \frac{1}{2}(X, Y)u\otimes[A, B]a = \\
  =[X, Y]u\otimes (A\circ B)a + \frac{1}{2}(XY + YX)u\otimes[A, B]a = XYu\otimes (A\circ B + \frac{1}{2}[A, B])a - \\
  -YXu\otimes (A\circ B - \frac{1}{2}[A, B])a =XYu\otimes ABa - YXu\otimes BAa =\\
  = [X\otimes A, Yu\otimes Ba] - [Y\otimes B, Xu\otimes Aa] = [R(\xi), R(\eta)]x.
\end{multline*}
Thus the map $R$ is a linear representation of the Lie algebra $\mathfrak{h}$. The truth of the Jacobi identity in case 2 follows from the defining relation of the representation $R$. It remains to prove the truth of the Jacobi identity for cases 3 and 4.

Consider case 3. In this case, the Jacobi identity is sufficient to prove for the following two sub-cases:
\begin{enumerate}
    \item for elements $\xi\in\mathfrak{g}_0, x, y\in\mathfrak{g}_1$;
    \item for elements $\xi\in\mathfrak{g}_2, x, y\in\mathfrak{g}_1$.
\end{enumerate}
In the first of the presented sub-cases of case 3, for $\xi=D\in\mathfrak{g}_0, x = u\otimes a, y=v\otimes b\in\mathfrak{g}_1$, the Jacobi identity looks like this:
\begin{equation*}
    \bigl[[D, u\otimes a], v\otimes b\bigr] + \bigl[[u\otimes a, v\otimes b], D\bigr] + \bigl[[v\otimes b, D], u\otimes a\bigr] = 0.
\end{equation*}
Applying commutational formulas (\ref{kom''1})--(\ref{kom''5}), the equation above can be rewritten as:
\begin{multline}\label{jacident}
S(u, v)\otimes\varphi(Da, b) + \langle u, v\rangle\delta(Da, b) -\\ - S(u, v)\otimes [D, \varphi(a, b)] - \langle u, v\rangle[D, \delta(a, b)] + \\ + S(u, v)\otimes\varphi(a, Db) + \langle u, v\rangle\delta(a, Db) = 0.
\end{multline}
Applying $\mathfrak{g}_0$-equivariance of $\delta$ for equation (\ref{jacident}), we get:
$$S(u, v)\otimes\varphi(Da, b)- S(u, v)\otimes [D, \varphi(a, b)]+ S(u, v)\otimes\varphi(a, Db) = 0.$$
Obviously, the resulting equation is equivalent to $\mathfrak{g}_0$-equivariancy of $\varphi$.

For further reasoning, we will need to derive two auxiliary identities. For arbitrary $A, B\in J_2, a, b\in J_1$ we have:
\begin{multline*}
    \bigl([A,\delta(a, b)], B\bigr) = \bigl([A,\delta_c(a, b)], B\bigr)= \\= \bigl(\delta_c(a, b), [B, A]\bigr) = \bigl\langle[B, A]a, b\bigl\rangle = \\ =\bigl\langle BAa,b\bigr\rangle - \bigl\langle ABa, b\bigr\rangle = \bigl\langle BAa, b\bigr\rangle - \bigl\langle Ba, Ab\bigr\rangle = \\ = \bigl(\varphi(Aa, b) - \varphi(a, Ab), B\bigr).
\end{multline*}
Since the equation above is fulfilled for an arbitrary $B\in J_2$, it follows from the non-degeneracy of scalar multiplication on $J_2$ that
\begin{equation}\label{id_1}
    [A, \delta(a, b)] = \varphi(Aa, b) - \varphi(a, Ab),\qquad\forall A\in J_2, a, b\in J_1.
\end{equation}
This is one of the identities we need. Now we derive the second identity. For arbitrary $A, B, C\in J_2, a, b\in J_1$ we have
\begin{multline*}
    \bigl([A,\varphi(a, b)], [B, C]\bigr) = \bigl(\varphi(a, b), [[B, C], A]\bigr) = \bigl\langle[[B, C], A]a, b\bigl\rangle = \\ =\bigl\langle [B, C]Aa,b\bigr\rangle - \bigl\langle A[B, C]a, b\bigr\rangle = \langle [B, C]Aa,b\rangle - \langle[B,C]a, Ab\rangle = \\ =\bigl(\delta_c(Aa, b) - \delta_c(a, Ab), [B, C]\bigr).
\end{multline*}
Since the equation above holds for all elements generating the Lie algebra $[J_2, J_2]$, then
\begin{equation*}\label{id_2}
    [A, \varphi(a, b)] = \delta_c(Aa, b) - \delta_c(a, Ab),\qquad\forall A\in J_2, a, b\in J_1.
\end{equation*}
Therefore, using the identity (\ref{delta0}) and the definition of $\delta$, we will get
\begin{equation}\label{id_2}
    [A, \varphi(a, b)] = \delta(Aa, b) - \delta(a, Ab),\qquad\forall A\in J_2, a, b\in J_1.
\end{equation}
This is the second of the identities we need.

Consider the Jacobi identity in the second of the above sub-cases of case 3 for $\xi = X\otimes A\in\mathfrak{g}_2, x = u\otimes a, y = v\otimes b\in\mathfrak{g}_1$:
\begin{equation*}
    \bigl[[X\otimes A, u\otimes a], v\otimes b\bigr] + \bigl[[u\otimes a, v\otimes b], X\otimes A\bigr] + \bigl[[v\otimes b, X\otimes A], u\otimes a\bigr] = 0.
\end{equation*}
Using the commutational relations (\ref{kom''1})--(\ref{kom''5}), the equation above can be rewritten in the form:
\begin{multline}\label{jacident2}
S(Xu, v)\otimes\varphi(Aa, b) + \langle Xu, v\rangle\delta(Aa, b) + \\ + [S(u, v), X]\otimes(\varphi(a, b)\circ A) + \frac{1}{2}(S(u, v), X)[\varphi(a, b), A] + \\ + \langle u, v\rangle X\otimes[\delta(a, b), A] - S(Xv, u)\otimes\varphi(Ab, a) - \langle Xv, u\rangle\delta(Ab, a) = 0.
\end{multline}
Using the identities (\ref{id_1}) and (\ref{id_2}), which have been proved above,  and lemma \ref{lem3}, the equation (\ref{jacident2}) can be rewritten in the form:
\begin{multline}\label{jid}
\bigl[X,S(u,v)\bigr]\otimes\left(\frac{1}{2}\varphi(Aa,b) + \frac{1}{2}\varphi(a,Ab) - \varphi(a,b)\circ A\right) = 0.
\end{multline}

To prove equation (\ref{jid}), we need to derive another auxiliary identity. Using the properties of scalar multiplication on the Jordan algebra $J_2$ and the equation $\varphi = \pi_2\varphi_m$, we have for arbitrary $A, B\in J_2, a, b\in J_1$:
\begin{multline*}
    (A\circ\varphi(a, b), B)  = (\varphi(a, b), A\circ B) =\\=\frac{1}{2}(AB, \varphi(a, b)) + \frac{1}{2}(BA,\varphi(a, b)) = \frac{1}{2}\langle ABa, b\rangle + \frac{1}{2}\langle BAa, b\rangle= \\= \frac{1}{2}\langle Ba, Ab\rangle + \frac{1}{2}\langle BAa, b\rangle= \frac{1}{2}(B, \varphi(a, Ab)) + \frac{1}{2}(B,\varphi(Aa, b)).
\end{multline*}
Then:
\begin{equation}\label{newid}
    A\circ\varphi(a, b) = \frac{1}{2}\bigl(\varphi(Aa, b) + \varphi(a, Ab)\bigr),\quad\forall A\in J_2, a, b\in J_1.
\end{equation}
Using the equation (\ref{newid}), we will get the identity (\ref{jid}) and the Jacobi identity in case 3.

Consider the case 4. For $x = u\otimes a, y = v\otimes b, z=w\otimes c\in\mathfrak{g}_1$ the Jacobi identity looks like this:
\begin{equation*}
    \bigl[[u\otimes a, v\otimes b], w\otimes c\bigr] + \bigl[[v\otimes b, w\otimes c], u\otimes a\bigr] + \bigl[[w\otimes c, u\otimes a], v\otimes b\bigr] = 0.
\end{equation*}
Using the commutational relations (\ref{kom''1})--(\ref{kom''5}), it can be rewritten in the following form:
\begin{equation}\label{jid4}
\langle u, v\rangle w\otimes F(a, b, c) + \langle v, w\rangle u\otimes F(b, c, a) + \langle w, u\rangle v\otimes F(c, a, b) = 0,
\end{equation}
where $$F(a, b, c) = \delta(a, b)c + \varphi(a, c)b + \varphi(b, c)a.$$

Since by definition of the symplectic Li-Jordan structure $(J_1; J_2; \mathfrak{g}_0; \delta_0)$ the map $F: J_1\times J_1\times J_1\rightarrow J_1$ is symmetric, then the equation (\ref{jid4}) can be rewritten as:
$$\bigl(\langle u, v\rangle w + \langle v, w\rangle u + \langle w, u\rangle v\bigr)\otimes F(a, b, c) = 0.$$
Using the lemma \ref{lem3} we may write $\langle u, v\rangle w + \langle v, w\rangle u + \langle w, u\rangle v = 0$, and then the Jacobi identity in the case 4 is fulfilled.

Thus, the vector space $\mathfrak{g}$ is endowed with a Lie algebra structure with a short $SL_2$-structure. We prove that the constructed Lie algebra is simple. Let $\mathfrak{j}$ be the ideal of the Lie algebra $\mathfrak{g}$. Obviously, $$\mathfrak{j} = \mathfrak{j}_0\oplus\mathfrak{j}_1\oplus\mathfrak{j}_2, \text{~where~}\mathfrak{j}_i =\mathfrak{g}_i\cap\mathfrak{j},\quad i = 0,1,2.$$

Consider $\mathfrak{j}_2$. It is clear that $\mathfrak{j}_2 = \mathfrak{sl}_2\otimes I_2$, with $I_2\triangleleft J_2$. Since the Jordan algebra $J_2$ is simple, two cases are possible:

1. $I_2 = 0$. Then $\mathfrak{j} = \mathfrak{j}_0\oplus\mathfrak{j}_1$. We may note, that $\mathfrak{j}_1 = \mathbb{C}^2\otimes I_1$, where $I_1\subset J_1$.

Since $\mathfrak{j}$ is invariant with respect to commuting with $\mathfrak{g}_1$, it follows from the identity (\ref{kom''3}) that $$\varphi(a,b) = 0,\qquad\forall a\in I_1, b\in J_1.$$
Therefore, $$\langle A a, b\rangle = (A, \varphi(a, b)) = 0,\quad\forall A\in J_2, a\in I_1, b\in J_1.$$
Using the non-degeneracy of scalar multiplication on $J_1$, it can be deduced from the identity above that $A a = 0,\,\forall A\in J_2, a\in I_1.$
But then $I_1 = 0$ due to the fact that the algebra $J_2$ contains an identical operator. So $\mathfrak{j} = \mathfrak{j}_0$.

Since $\mathfrak{j}$ is invariant with respect to commuting with $\mathfrak{g}_1$, then $Da = 0,\,\forall D\in\mathfrak{j}_0, a\in J_1$, from which it immediately follows that $\mathfrak{j} = 0$.

2. $I_2 = J_2$. Then $\mathfrak{j} = \mathfrak{j}_0\oplus\mathfrak{j}_1\oplus\mathfrak{g}_2$. Since $\mathfrak{j}$ is invariant with respect to commuting with $\mathfrak{g}_1$, and the algebra $J_2$ contains an identical operator, then $\mathfrak{j}_1 = \mathfrak{g}_1$. Since the map $\delta$ is surjective, $\mathfrak{g}_0$ is generated by the operators $\delta(a,b)$, where $a, b\in J_1$. Hence, since the ideal $\mathfrak{j}$ contains all such operators, then $\mathfrak{j}_0 = \mathfrak{g}_0$ and $\mathfrak{j} =\mathfrak{g}$.

Thus, the constructed Lie algebra $\mathfrak{g}$ is simple, so the theorem is proved.
\end{Proof}

\section{Classification of short $SL_2$-structures}\label{classif}
\subsection{Preliminaries}\label{premil}

Consider an arbitrary short $SL_2$-structure on a simple Lie algebra $\mathfrak{g}$ and keep all the notation of the previous sections, and through $e,f, h$ we again denote the basic elements of the Lie algebra $\mathfrak{sl}_2\subset\mathfrak{g}$ satisfying the relations (\ref{sl2}).

Note that due to the fact that $\mathfrak{g}_1 = \mathbb{C}^2\otimes J_1\neq 0$, the only possible nonzero value of the simple root of the Lie algebra $\mathfrak{g}$ on the element $h$ is 1. In the Dynkin diagram of the Lie algebra $\mathfrak{g}$, we will black out the vertices corresponding to simple roots having the value 1 on the element $h$, and we will call such vertices {\it unit vertices}. The vertices of the diagram corresponding to the roots having a zero value on the $h$ element will be drawn uncolored and called {\it zero vertices}. The resulting diagram is called {\it diagram of a short $SL_2$-structure on the Lie algebra $\mathfrak{g}$}.

Using the diagram of a short $SL_2$-structure, it is easy to calculate the dimension of the center $\mathfrak{z}^0$ of the algebra $\mathfrak{g}^0$. Since the Cartan subalgebra is $\mathfrak{t}\subset\mathfrak{g}^0$, the dimension of the center is equal to the difference between the rank of the Lie algebra $\mathfrak{g}$ and the rank of the semisimple part of the algebra $\mathfrak{g}^0$, that is, the number of black vertices.

The \ref{finth} theorem allows us to assert that by enumerating all existing short $SL_2$-structures on simple Lie algebras and specifying simple Lie-Jordan symplectic structures corresponding to these structures, we thereby enumerate all possible simple Lie-Jordan symplectic structures, which, generally speaking, can be considered by themselves, without binding them to $SL_2$-structures. In order to do this, we denote by $\widetilde{G}^0$ an algebraic group whose tangent algebra is the Lie algebra $\tilde{\mathfrak{g}}^0$, and use the following statement:
\begin{theorem}\label{invth}
For a given $\mathbb{Z}$-grading defined by a semisimple element $h\in\mathfrak{g}$, the following conditions are equivalent:
\begin{enumerate}
\item The element $h$  is included into the $\mathfrak{sl}_2$-triple;
\item The representation $\widetilde{G}^0:\mathfrak{g}^2$ has no open orbits;
\item On the subspace $\mathfrak{g}^2$ there is a nontrivial polynomial invariant of the action $\widetilde{G}^0:\mathfrak{g}^2$.
\end{enumerate}
\end{theorem}
The proof of equivalence of points 1 and 2 from the theorem above can be found in \cite[point~1,~section~1.4]{Vinberg_1}. Point 2 obviously follows from point 3. The proof that point 1 follows point 3 can be found in \cite[paragraph~1,~propositions~1.1~and~1.2]{Kac}.

It is important to note that $\mathfrak{sl}_2$ is a triple containing a semisimple element $h$, defined uniquely, up to conjugation with $G^0$, where $G^0$ is a connected algebraic group with tangent algebra $\mathfrak{g}^0$. Namely, as an element $e$, we can take any element from the open $G^0$-orbit in $\mathfrak{g}^2$, and the element $f$ is uniquely determined by the elements $h$ and $e$.

Now we have everything we need to classify all short $SL_2$-structures on simple Lie algebras.

\subsection{Classification of short $SL_2$-structures on classical simple Lie algebras}\label{classalg}

In the case of classical simple Lie algebra, the element $h$ is a diagonal matrix, and diagonal elements can be calculated from the values of $h$ on simple roots. By commuting $h$ with an arbitrary matrix from the algebra $\mathfrak{g}$, one can define the representation space $\mathfrak{g}^2$, as well as the group $\widetilde{G}^0$. This method completely determines the action of $\widetilde{G}^0:\mathfrak{g}^2$, and using the theorem  \ref{invth}, or rather its third point, it is possible to determine whether there is a short $SL_2$-structure in this case and, if there is, under what conditions it exists. Next, we will denote everywhere the highest root of the algebra $\mathfrak{g}$ by $\alpha_{\operatorname{max}}$, and simple roots by $\alpha_i$.

{\bf The case $\bf\mathfrak{g} = \mathfrak{sl}_{n}$}. The dependence of the highest root from simple roots is expressed by the formula $$\alpha_{\operatorname{max}} = \alpha_1 + ... + \alpha_{n-1}.$$ Since $\alpha_{\operatorname{max}}(h) = 2$ and the set of single vertices of the diagram is not empty, then a single case is possible:
\begin{equation*}
\alpha_i(h) = 1, \alpha_j(h) = 1, \alpha_k(h) = 0,\qquad\forall k = \overline{1, (n-1)}, k\neq i, j, i<j.
\end{equation*}
In this case, the element $h$ is a diagonal matrix of size $n\times n$ of the following form (then everywhere through $E^{(i)}$ we will denote a unit matrix of size $i\times i$):

$$h =  \begin{pmatrix}
cE^{(i)} &   0   & 0 \\
0  & (c-1)E^{(j-i)} & 0 \\
0  &    0   & (c-2)E^{(n-j)}
\end{pmatrix}.$$
The multiplier $c\in\mathbb{C}$ in the formula above is uniquely determined from the condition $\operatorname{tr}(h) = 0$.

The algebra $\mathfrak{g}$ can be conditionally represented as:

$$\mathfrak{g} =  \begin{pmatrix}
\mathfrak{g}^{0}_{(1)} &   \mathfrak{g}^{1}_{(1)}    & \mathfrak{g}^{2} \\
\mathfrak{g}^{-1}_{(1)}  & \mathfrak{g}^{0}_{(2)} & \mathfrak{g}^{1}_{(2)} \\
\mathfrak{g}^{-2}  &    \mathfrak{g}^{-1}_{(2)}   & \mathfrak{g}^{0}_{(3)}
\end{pmatrix}.$$

This formula means that the algebra $\mathfrak{g}$ is represented as a direct sum of subspaces $\mathfrak{g}^{k}$ for $k\neq 0$, which are spaces of matrices from $\mathfrak{g}$, in which nonzero elements can only stand in places, corresponding to the location of the blocks $\mathfrak{g}^k_{(m)}$ for $m = 1, 2$ in the formula above, and the subspace $\mathfrak{g}^0$, which is a space of block-diagonal matrices from $\mathfrak{sl}_n$, consisting of blocks $\mathfrak{g}^0_{(m)}$ for $m = 1, 2, 3$. Here the block $\mathfrak{g}^{0}_{(1)}$ has the size $i\times i$, block $\mathfrak{g}^{0}_{(2)}$ has the size $(j-i)\times(j-i)$ and block $\mathfrak{g}^{0}_{(3)}$ has the size $(n-j)\times(n-j)$. The superscript of each block indicates its own space to which this block belongs. The subscripts are designed to distinguish between different blocks belonging to the same proper space. Let 's draw a diagram of this case:
\begin{center}
\begin{dynkinDiagram} [labels*= {,,\alpha_{i},,,\alpha_{j}},root radius=.10cm, edge length=1.3cm]A{o.o*o.o*o.o}
\end{dynkinDiagram}
\end{center}
Note that the semisimple part of $\mathfrak{g}^0$ coincides with $\mathfrak{sl}_i\oplus\mathfrak{sl}_{j-i}\oplus\mathfrak{sl}_{n-j}$, as well as $\dim\mathfrak{z}^0 = 2$, whence we conclude that $\widetilde{G}^{0} = SL_i\times SL_{j-i}\times SL_{n-j}\times Z$, where $Z$ is a one-dimensional center. It is also easy to conclude that $\mathfrak{g}^2 = \mathbb{C}^i\otimes(\mathbb{C}^{n-j})^{*}$ as $\widetilde{G}^0$-module. Since the action of the one-dimensional center $Z$ on $\mathfrak{g}^2$ is trivial, it is obvious that an invariant on $\mathfrak{g}^2$ for the group $\widetilde{G}^{0}$ exists if and only if $i=n-j$, and is equal to the determinant of the matrix from $\mathfrak{g}^2$.

We define a simple symplectic Lie-Jordan structure that corresponds to this structure. It is easy to understand that the elements of $e,f,h$ for a given short $SL_2$-structure are matrices of size $n\times n$ of the following form:
$$h =  \begin{pmatrix}
E^{(i)} &   0   & 0 \\
0  & 0 & 0 \\
0  &    0   & -E^{(i)}
\end{pmatrix},\quad e =  \begin{pmatrix}
0 &   0   & E^{(i)} \\
0  & 0 & 0\\
0  &    0   & 0
\end{pmatrix},\quad f =  \begin{pmatrix}
0 &   0   & 0 \\
0  & 0 & 0 \\
E^{(i)}  &    0   & 0
\end{pmatrix}.$$
It is obvious, that $$\mathfrak{g}_0 = \mathfrak{z}_{\mathfrak{g}^0}(\mathfrak{sl}_2), \mathfrak{g}^1 = e_1\otimes J_1, \mathfrak{g}^2 = e\otimes J_2,$$
where $$\mathfrak{z}_{\mathfrak{g}^0}(\mathfrak{sl}_2) = \{\xi\in\mathfrak{g}^{0}: [\xi, \eta] = 0,\quad\forall\eta\in\mathfrak{sl}_2\subset\mathfrak{g}_2\}.$$
Let's denote
$$\mathfrak{c}_{i,n} = \{T\in\mathfrak{gl}_n: T=\operatorname{diag}\{\underbrace{c_1,..,c_1}_{i},\underbrace{c_2,..,c_2}_{n-2i}, \underbrace{c_1,..,c_1}_{i}\},\,2c_1i + c_2(n-2i) = 0,\,c_1,c_2\in\mathbb{C}\}.$$
From the equations above, it is not difficult to conclude that:
$$J_1 = \mathfrak{g}^1_{(1)}\oplus\mathfrak{g}^1_{(2)} \cong \left(\mathbb{C}^{i}\otimes\bigl(\mathbb{C}^{n-2i}\bigr)^{*}\right)\oplus\left(\mathbb{C}^{n-2i}\otimes\bigl(\mathbb{C}^{i}\right)^{*}\bigr),$$
\begin{equation*}
  \mathfrak{g}^0_{(1)}\oplus\mathfrak{g}^0_{(2)}\oplus\mathfrak{g}^0_{(3)}\supset\mathfrak{g}_0 = \{(A, B, A), \text{~where~} A\in\mathfrak{gl}_{i}, B\in\mathfrak{gl}_{n-2i},\,2\operatorname{tr}A + \operatorname{tr}B = 0\},
\end{equation*}
$$\mathfrak{g}_0\cong\mathfrak{sl}_i\oplus\mathfrak{sl}_{n-2i}\oplus\mathfrak{c}_{i,n},$$
$$J_2 =\mathfrak{g}^{2}\cong\mathfrak{gl}_i.$$
Let's define scalar multiplication by $\mathfrak{g}=\mathfrak{sl}_n$ using the following formula:
\begin{equation*}
  (A, B) = (n - 2i)\operatorname{tr}(AB),\quad\forall A,B\in\mathfrak{sl}_n.
\end{equation*}
The normalization factor $n-2i$ in the formula above is chosen in such a way that, as a result of calculations, scalar multiplication on the Jordan algebra $J_2\cong\mathfrak{gl}_i$ would equal the trace of the product of operators on the space $J_1$. Let's check that this is really the case. The actions of the Lie algebra $\mathfrak{g}_0$ and the Jordan algebra $J_2 = \mathfrak{gl}_i$ on the space $J_1$ are carried out by the following formulas:
$$C(a, b) = (Ca, bC),\quad\forall C\in\mathfrak{gl}_i, (a, b)\in J_1;$$
$$(D_1, D_2) (a, b) = (D_1a - aD_2, D_2b-bD_1),\quad\forall (D_1, D_2)\in\mathfrak{g}_0, (a, b)\in J_1.$$
Therefore, the trace of the product of operators from $J_2$ has the form
\begin{equation*}
  \operatorname{tr}(\mathcal{A}_1\mathcal{B}_1) = 2(n-2i)\operatorname{tr}(AB),\quad\forall A,B\in\mathfrak{gl}_i,
\end{equation*}
where $\mathcal{C}_1$ denotes a linear operator on the space $J_1$ corresponding to the matrix $C\in\mathfrak{gl}_i$.

Considering an arbitrary element $\tilde{e}\otimes A\in\mathfrak{g}^2 = \tilde{e}\otimes J_2$ and scalar multiplying it by an arbitrary element $\tilde{f}\otimes B\in\mathfrak{g}^{-2} = \tilde{f}\otimes J_2$, we get that
$$(A, B)_2 = 2(n-2i)\operatorname{tr}(AB),\quad\forall A, B\in\mathfrak{gl}_i,$$
which confirms the above agreement.

With this choice of the normalization factor of scalar multiplication by $\mathfrak{g} = \mathfrak{sl}_n$, the symplectic structure on the space $J_1$ is given using skew-scalar multiplication, which has the form $\forall a_1, a_2\in\mathbb{C}^{i}\otimes\left(\mathbb{C}^{n-2i}\right)^{*}$,$b_1,b_2\in\mathbb{C}^{n-2i}\otimes\left(\mathbb{C}^{i}\right)^{*}$:
$$\bigl((a_1, b_1),(a_2, b_2)\bigr) = (n-2i)\operatorname{tr}(b_1a_2 - a_1b_2).$$
The multiplication operation on the Jordan algebra $J_2$ looks like:
$$A\circ B = \dfrac{1}{2}(AB + BA),\quad\forall A, B\in\mathfrak{gl}_i.$$
The action of the Lie algebra $\mathfrak{g}_0$ on the Jordan algebra $J_2$ is carried out by the following formula:
$$(D_1, D_2)(A) =  [D_1, A],\quad\forall (D_1, D_2)\in\mathfrak{g}_0, A\in J_2.$$
From here it is not difficult to get that
$$\mathfrak{i}_0 = \mathfrak{sl}_{n-2i}\oplus\mathfrak{c}_{i,n}, [J_2, J_2] = \mathfrak{sl}_i.$$

The maps $\varphi$ and $\delta$ in this case have the following form$\forall a_1, a_2\in\mathbb{C}^{i}\otimes\left(\mathbb{C}^{n-2i}\right)^{*},$ $b_1,b_2\in\mathbb{C}^{n-2i}\otimes\left(\mathbb{C}^{i}\right)^{*}$:
$$\varphi\bigl((a_1, b_1), (a_2, b_2)\bigr) = \dfrac{1}{2}(a_2b_1 - a_1b_2),$$
$$\delta\bigl((a_1, b_1), (a_2, b_2)\bigr) = \left(-\dfrac{1}{2}(a_1b_2 + a_2b_1), b_1a_2 + b_2a_1\right).$$
Also
$$\delta_0\bigl((a_1, b_1), (a_2, b_2)\bigr) = \left(-\dfrac{1}{2i}\operatorname{tr}(a_1b_2 + a_2b_1)E^{(i)}, b_1a_2 + b_2a_1\right).$$
Analyzing further cases, we will not describe in such detail all the actions and maps that occur in one or another short $SL_2$-structure. If desired, all of this can be done by direct calculations.

{\bf The case $\bf\mathfrak{g} = \mathfrak{so}_{2n+1}$}. We will consider the algebra $\mathfrak{g} = \mathfrak{so}_{2n+1}$ in a basis in which its elements are matrices of size $(2n+1)\times(2n+1)$, skew-symmetric with respect to the side diagonal. In this case, $$\alpha_{\operatorname{max}} = \alpha_1 + 2\alpha_2 + ... + 2\alpha_n.$$ Since the set of unit vertices in the diagram of a short $SL_2$-structure is not empty, the only possible case for this Lie algebra is the following case:
\begin{equation*}
\alpha_i(h) = 1,\qquad i=\overline{2, n}.
\end{equation*}
In this case, the element $h$ is the matrix $(2n+1)\times(2n+1)$ of the following form:
$$h = \begin{pmatrix}
E^{(i)}&0&0\\
0&0&0\\
0&0&-E^{(i)}
\end{pmatrix}$$

Similarly to the case of $\mathfrak{g} = \mathfrak{sl}_n$, the algebra $\mathfrak{g}= \mathfrak{so}_{2n+1}$ can be conditionally represented as:
$$\mathfrak{g} =  \begin{pmatrix}
\mathfrak{g}^{0}_{(1)} &   \mathfrak{g}^{1}    & \mathfrak{g}^{2} \\
\mathfrak{g}^{-1}  & \mathfrak{g}^{0}_{(2)} & -\bigl(\mathfrak{g}^{1}\bigr)^{\mathsf{s}} \\
\mathfrak{g}^{-2}  &    -\bigl(\mathfrak{g}^{-1}\bigr)^{\mathsf{s}}   & -\bigl(\mathfrak{g}^{0}_{(1)}\bigr)^{\mathsf{s}}
\end{pmatrix}.$$
Here is the block $\mathfrak{g}^{0}_{(1)}$ has the size $i\times i$, block $\mathfrak{g}^{0}_{(2)}$ has the size $(2(n-i)+1)\times(2(n-i)+1)$, and the remaining blocks have the appropriate sizes. The only difference from the previous case is that the terms in the decomposition of the algebra $\mathfrak{g}$ into a direct sum using the formula above in this case are: the space of block-diagonal matrices consisting of blocks $\mathfrak{g}^0_{(m)}$, the space of matrices, which nonzero elements can only stand in places inside the block $\mathfrak{g}^{1}$ and a block symmetric to it with respect to the side diagonal, in which there is a matrix antisymmetric to the matrix from the block $\mathfrak{g}^1$ with respect to the side diagonal, a similar summand for the block $\mathfrak{g}^{-1}$, and also two matrix spaces where nonzero elements can only stand in place inside the blocks $\mathfrak{g}^{2}$ and $\mathfrak{g}^{-2}$.
From this it is immediately clear that $\mathfrak{g}^2$ is a space of matrices of size $i\times i$, skew-symmetric with respect to the side diagonal, that is, we can conditionally write that $\mathfrak{g}^2 = \bigwedge^{2}\mathbb{C}^{i}$. Now we get $\widetilde{G}^{0}$. Consider the diagram:
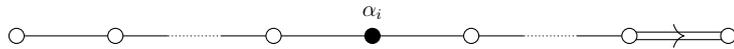
\begin{center}
\begin{dynkinDiagram} [labels*= {,,,\alpha_{i}},root radius=.10cm, edge length=1.3cm]B{oo.o*o.oo}
\end{dynkinDiagram}
\end{center}
It gives us on the one hand that $\dim\mathfrak{z}^0 =1$ and $\mathfrak{z}^0 = \langle h\rangle$, and on the other hand that the semisimple part of $\mathfrak{g}^0$ is equal to $\mathfrak{sl}_i\oplus\mathfrak{so}_{2(n-i)+1}$, so we can conclude that $\widetilde{G}^{0}= SL_{i}\times SO_{2(n-i)+1}$, and only the first multiplier acts on $\mathfrak{g}^2$ . The invariant of this action is the determinant of the matrix from $\mathfrak{g}^{2}$, which is not equal to 0 if and only if $i$ is even. Thus, a short $SL_2$-structure is possible only with an even $i$.

We define a simple symplectic Lie-Jordan structure that corresponds to a given $SL_2$-structure. Let's define scalar multiplication on $\mathfrak{g}=\mathfrak{so}_{2n+1}$ using the following formula:
\begin{equation*}
  (A, B) = \frac{1}{2}(2(n-i)+1)\operatorname{tr}(AB),\quad\forall A,B\in\mathfrak{so}_{2n+1}.
\end{equation*}
Since $i = 2k$, then $\mathfrak{sl}_2$-triple $e, f, h$ for this structure are matrices of size $(2n+1)\times (2n+1)$ of the following form:
$$h =  \begin{pmatrix}
E^{(i)} &   0   & 0 \\
0  & 0 & 0 \\
0  &    0   & -E^{(i)}
\end{pmatrix},\quad e =  \begin{pmatrix}
0 &   0   & I_{i} \\
0  & 0 & 0\\
0  &    0   & 0
\end{pmatrix},\quad f =  \begin{pmatrix}
0 &   0   & 0 \\
0  & 0 & 0 \\
I_{i}  &    0   & 0
\end{pmatrix},$$
where $I_{i}$ is a matrix of size $i\times i$ of the following form:
$$I_{i} =  \begin{pmatrix}
E^{(k)} & 0\\
0  & -E^{(k)}
\end{pmatrix}.$$
Acting similarly to the case of $\mathfrak{g} = \mathfrak{sl}_n$, it is not difficult to obtain that in this case:
$$\mathfrak{g}_0 = \mathfrak{sp}_i\oplus\mathfrak{so}_{2(n-i)+1},\quad J_1 = \mathbb{C}^{i}\otimes\left(\mathbb{C}^{2(n-i)+1}\right)^{*},\quad J_2=\mathfrak{sym}^{-}_i,\qquad i=2k.$$
Here, using $\mathfrak{sym}^{-}_i$ for $i=2k$, the Jordan algebra of matrices of order $i$ symmetric with respect to a skew-symmetric non-degenerate bilinear form is denoted.
Denote
$$a = \begin{pmatrix}
a_1\\
a_2
\end{pmatrix}\in J_1,\quad b = \begin{pmatrix}
b_1\\
b_2
\end{pmatrix}\in J_1,$$
where $a_i, b_i$ are arbitrary matrices of size $k\times (2(n-i)+1)$.
The symplectic structure on the space $J_1$ is given by skew-scalar multiplication, which has the form:
$$\bigl(a,b\bigr) = (2(n-i)+1)\operatorname{tr}(a_1b_2^{\mathsf{s}} - a_2b_1^{\mathsf{s}}).$$
The multiplication operation on the Jordan algebra $J_2$ looks like:
$$A\circ B = \dfrac{1}{2}(AB + BA),\quad\forall A, B\in\mathfrak{sym}^{-}_i.$$
The action of the Jordan algebra $J_2 = \mathfrak{sym}^{-}_i$ on the space $J_1$ is carried out by multiplying the matrix from $J_1$ on the left by the corresponding matrix from $\mathfrak{g}^2$. The action of the Lie algebra $\mathfrak{g}_0$ on the space $J_1$ is carried out by the formula:
$$(D_1, D_2)a = D_1a - aD_2,\quad\forall D_1\in\mathfrak{sp}_i, D_2\in\mathfrak{so}_{2(n-i) + 1}, a\in J_1.$$

{\bf The case $\bf\mathfrak{g} = \mathfrak{so}_{2n}$}.We will consider the algebra $\mathfrak{g} = \mathfrak{so}_{2n}$ in a basis in which its elements are matrices of size $(2n)\times(2n)$, skew-symmetric with respect to the side diagonal. Scalar multiplication on $\mathfrak{so}_{2n}$ will be set similarly to the previous case. For this algebra, the dependence of the highest root on the simples looks like this: $$\alpha_{\operatorname{max}} = \alpha_1 + 2\alpha_2 + ... + 2\alpha_{n-2} + \alpha_{n-1} + \alpha_{n}.$$
Here, generally speaking, $4$ cases are possible:
\begin{enumerate}
\item $\alpha_i(h) = 1, \alpha_j(h) = 0,\qquad j = \overline{1,n}, i=\overline{2,(n-2)}, j\neq i$;
\item $\alpha_n(h) = \alpha_{n-1}(h) = 1, \alpha_i(h) = 0,\qquad i = \overline{1,n}, i\neq n, n-1$;
\item $\alpha_1(h) = \alpha_n(h) = 1, \alpha_i(h) = 0,\qquad i = \overline{1,n}, i\neq 1, n$;
\item $\alpha_1(h) = \alpha_{n-1}(h) = 1, \alpha_i(h) = 0,\qquad i = \overline{1,n}, i\neq 1, n-1$.
\end{enumerate}

Cases 3 and 4 are translated into each other using the automorphism of the Dynkin diagram, so it is sufficient to consider only one of them. Thus, it remains to consider cases 1, 2 and 3.

In case 1, the matrix $h$ will have the same form as in the case of $\mathfrak{g} = \mathfrak{so}_{2n+1}$ with the only difference that the central zero matrix here has dimensions $2(n-i)\times2(n-i)$. Therefore, the answer here will be the same: an invariant (determinant) exists if and only if $i$ is even (here $\mathfrak{g}^2 =\bigwedge^{2}\mathbb{C}^i, \widetilde{G}^0= SL_i\times SO_{2(n-i)}$). The diagram of this case looks like this:
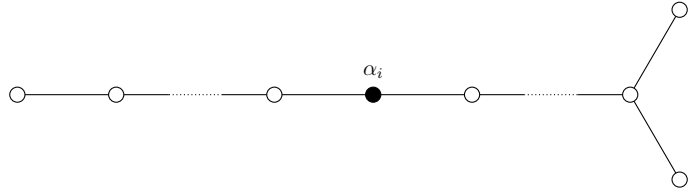
\begin{center}
\begin{dynkinDiagram} [labels*= {,,,\alpha_{i}},root radius=.10cm, edge length=1.3cm]D{oo.o*o.ooo}
\end{dynkinDiagram}
\end{center}
Here $\mathfrak{sl}_2$-the triple $e, f, h$ are matrices of size $2n\times 2n$ similar to the case of $\mathfrak{g} =\mathfrak{so}_{2n+1}$ of the form. Therefore, acting in the same way as for $\mathfrak{so}_{2n+1}$, we have:
$$\mathfrak{g}_0 = \mathfrak{sp}_i\oplus\mathfrak{so}_{2(n-i)},\quad J_1 = \mathbb{C}^{i}\otimes\left(\mathbb{C}^{2(n-i)}\right)^{*},\quad J_2=\mathfrak{sym}^{-}_i,\qquad i=2k\leqslant n-2.$$
Skew-symmetric scalar multiplication on the space $J_1$, as well as the action of the Lie algebra $\mathfrak{g}_0$ and the Jordan algebra $J_2$ on the space $J_1$ in this case are carried out by formulas similar to the case $\mathfrak{g} = \mathfrak{so}_{2n+1}$.

Consider case 2. In this case $$h = \operatorname{diag}\{1,...,1,0,0,-1,...,-1\}.$$
Algebra $\mathfrak{g} = \mathfrak{so}_{2n}$ similarly to the previous case can be conditionally represented as:
$$\mathfrak{g} = \begin{pmatrix}
\mathfrak{g}^{0}_{(1)} & \mathfrak{g}^{1}_{(1)} & \mathfrak{g}^{2} \\
\mathfrak{g}^{-1}_{(1)} & \mathfrak{g}^{0}_{(2)} & -\bigl(\mathfrak{g}^{1}_{(1)}\bigr)^{\mathsf{s}} \\
\mathfrak{g}^{-2} & -\bigl(\mathfrak{g}^{-1}_{(1)}\bigr)^{\mathsf{s}} & -\bigl(\mathfrak{g}^{0}_{(1)}\bigr)^{\mathsf{s}}
\end{pmatrix}.$$
Here the block $\mathfrak{g}^{0}_{(1)}$ has the size $(n-1)\times(n-1)$, block $\mathfrak{g}^{0}_{(2)}$ has the size $2\times 2$, and the remaining blocks have the corresponding sizes.

From the formula above it can be seen that $\mathfrak{g}^2 = \bigwedge^2\mathbb{C}^{n-1}$. Consider the diagram of this representation:
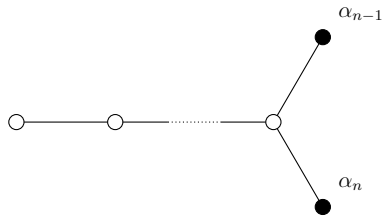
\begin{center}
\begin{dynkinDiagram} [labels*= {,,,\alpha_{n-1},\alpha_n},root radius=.10cm, edge length=1.3cm]D{oo.o**}
\end{dynkinDiagram}
\end{center}

From the diagram we deduce that the semisimple part of $\mathfrak{g}^0$ is equal to $\mathfrak{sl}_{n-1}$ and $\dim\mathfrak{z}^0 = 2$. It is not difficult to make sure that the vector orthogonal to the vector $h$ in $\mathfrak{z}^0$ corresponds to the block $\mathfrak{g}^{0}_{(2)}$ and, due to the locations of the blocks $\mathfrak{g}^{0}_{(2)}$ and $\mathfrak{g}^2$, acts on $\mathfrak{g}^2$ trivially, we can assume that $\widetilde{G}^{0}= SL_{n-1}$. Thus, it is obvious that this representation has an invariant (equal to the determinant of the matrix from $\mathfrak{g}^2$) if and only if $n$ is odd. That is, for an odd $n$ in this case, there is a short $SL_2$-structure.

Here the matrices $e, f, h$ look similar to the previous case. Therefore:
$$\mathfrak{g}_0 = \mathfrak{sp}_{n-1}\oplus\mathfrak{so}_{2},\quad J_1=\mathbb{C}^{n-1}\otimes\left(\mathbb{C}^{2}\right)^{*},\quad J_2=\mathfrak{sym}^{-}_{n-1},\qquad n=2k+1.$$
Skew-symmetric scalar multiplication on the space $J_1$, as well as the action of the Lie algebra $\mathfrak{g}_0$ and the Jordan algebra $J_2$ on the space $J_1$ in this case are carried out by formulas similar to the case $\mathfrak{g} = \mathfrak{so}_{2n+1}$.

Let's move on to the last case 3. The element $h$  in this case has the form $$h = \operatorname{diag}\left\{\frac{3}{2},\frac{1}{2},...,\frac{1}{2},-\frac{1}{2},...,-\frac{1}{2},-\frac{3}{2}\right\}.$$
Similarly to the previous cases, the algebra $\mathfrak{g} = \mathfrak{so}_{2n}$ can be conditionally represented as:
$$\mathfrak{g} = \begin{pmatrix}
\mathfrak{g}^{0}_{(1)} & \mathfrak{g}^{1}_{(1)} & \mathfrak{g}^{2} & 0\\
\mathfrak{g}^{-1}_{(1)} & \mathfrak{g}^{0}_{(2)} & \mathfrak{g}^{1}_{(2)} & -\bigl(\mathfrak{g}^{2}\bigr)^{\mathsf{s}} \\
\mathfrak{g}^{-2} & \mathfrak{g}^{-1}_{(2)} & -\bigl(\mathfrak{g}^{0}_{(2)}\bigr)^{\mathsf{s}} & -\bigl(\mathfrak{g}^{1}_{(1)}\bigr)^{\mathsf{s}} \\
0 & -\bigl(\mathfrak{g}^{-2}\bigr)^{\mathsf{s}} & -\bigl(\mathfrak{g}^{-1}_{(1)}\bigr)^{\mathsf{s}} &-\bigl(\mathfrak{g}^{0}_{(1)}\bigr)^{\mathsf{s}}
\end{pmatrix}.$$
Here $\mathfrak{g}^{0}_{(2)}, \mathfrak{g}^{1}_{(2)}, \mathfrak{g}^{-1}_{(2)}$ are blocks of size $(n-1)\times(n-1)$, $\mathfrak{g}^{1}_{(1)}, \mathfrak{g}^{2}$ are strings of length $n-1$, and $\mathfrak{g}^{-2}, \mathfrak{g}^{-1}_{(1)}$ are columns of height $n-1$. From where we get that $\mathfrak{g}^2 = \mathbb{C}^{n-1}$.

Consider the diagram of this case:
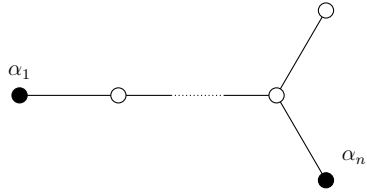
\begin{center}
\begin{dynkinDiagram} [labels*= {\alpha_1,,,,\alpha_{n}},root radius=.10cm, edge length=1.3cm]D{*o.oo*}
\end{dynkinDiagram}
\end{center}
The diagram shows that $\dim\mathfrak{z}^0 =2$, and the semisimple part of $\mathfrak{g}^0$ is equal to $\mathfrak{sl}_{n-1}$. Then $\widetilde{G}^0 = GL_{n-1}$ and we get a tautological representation of the group $GL_{n-1}$ on $(n-1)$-dimensional vector space, which obviously has no nontrivial invariants, that is, in this case does not exist short $SL_2$-structures.

{\bf The case $\bf\mathfrak{g} = \mathfrak{sp}_{2n}$}. We will consider the algebra $\mathfrak{sp}_{2n}$ in a basis in which its elements are given by matrices of the form:
$$\begin{pmatrix}
X&Y\\
Z&-X^{\mathsf{s}}
\end{pmatrix}\in\mathfrak{gl}_{2n}: X, Y, Z\in\mathfrak{gl}_n,\,Y^{\mathsf{s}} = Y,\,Z^{\mathsf{s}} = Z.$$
Scalar multiplication on $\mathfrak{sp}_{2n}$ will be set similarly to the previous case.

For this algebra, the highest root has the form $$\alpha_{\operatorname{max}} = 2\alpha_1 + ... + 2\alpha_{n-1} + \alpha_{n},$$ thus, a single case is possible here:
$$
\alpha_i(h) = 1, \alpha_j(h) = 0,\qquad i=\overline{1,(n-1)}, j =\overline{1,n}, j\neq i.
$$
The elements $e, f, h$ are matrices of size $2n\times 2n$ of the form similar to the case of $\mathfrak{g} = \mathfrak{sl}_n$, so when commuting, the subspace $\mathfrak{g}^2$ will correspond to the upper right corner block of size $i\times i$, from where it is clear that $\mathfrak{g}^2$ is the space of matrices of size $i\times i$ symmetric with respect to the side diagonal. Consider the diagram of this case:
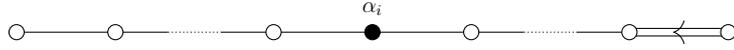
\begin{center}
\begin{dynkinDiagram} [labels*= {,,,\alpha_{i}},root radius=.10cm, edge length=1.3cm]C{oo.o*o.oo}
\end{dynkinDiagram}
\end{center}
Based on it, we get that $\dim\mathfrak{z}^0 =1$, and the semisimple part of $\mathfrak{g}^0$ coincides with $\mathfrak{sl}_i\oplus\mathfrak{sp}_{2(n-i)}$, from where $\widetilde{G}^0 = SL_i\times Sp_{2(n-i)}$. Then, since $Sp_{2(n-i)}$ acts on $\mathfrak{g}^2$ trivially, we can assume that $\widetilde{G}^0 = SL_i$, from which it is clear that in this case the determinant of the matrix of $\mathfrak{g}^2$ is the invariant of the action, that is, in this case, a short $SL_2$-structure exists for any $i=\overline{1,(n-1)}$.
In this case:
$$\mathfrak{g}_0 = \mathfrak{so}_{i}\oplus \mathfrak{sp}_{2(n-i)},\quad J_1 = \mathbb{C}^{i}\otimes\left(\mathbb{C}^{2(n-i)}\right)^{*},\quad J_2=\mathfrak{sym}^{+}_{i},\qquad i<n.$$
Here, $\mathfrak{sym}^{+}_i$ denotes the Jordan algebra of matrices of order $i$ symmetric with respect to a symmetric non-degenerate bilinear form.

Skew-symmetric scalar multiplication on the space $J_1$, as well as the action of the Lie algebra $\mathfrak{g}_0$ and the Jordan algebra $J_2$ on the space $J_1$ in this case are also carried out by formulas similar to the case $\mathfrak{g} = \mathfrak{so}_{2n+1}$.
\subsection{Classification of short $SL_2$-structures on special simple Lie algebras}\label{classspec}

Although each of the special simple Lie algebras has its own matrix model, all these models are quite cumbersome, so the study of short $SL_2$-structures on them in the way described above will be too complicated.

Consider the irreducible representation of the Lie algebra $\tilde{\mathfrak{g}}^{0}$ on $\mathfrak{g}^2$. It is uniquely determined by its highest weight, namely its decomposition by fundamental weights. The highest weight of this representation obviously coincides with the highest root of the Lie algebra $\mathfrak{g}$. Its decomposition by fundamental weights is obtained by calculating Cartan products between the highest root of the Lie algebra $\mathfrak{g}$ and the simple roots $\mathfrak{g}$ corresponding to the zero vertices of the diagram of the short $SL_2$-structure. Therefore, using only the extended Dynkin diagram, it is possible to define both the group $\widetilde{G}^0$ and the representation space $\mathfrak{g}^2$.

For further reasoning, we need the following notation. We will denote by $S^k$ the space of the spinor representation of the algebra $\mathfrak{so}_k$, and by $S^k_{1/2}$ and $S^k_{-1/2}$ --- the spaces of the semi-spinor representations of the corresponding algebra for even $k$. When the algebra $\mathfrak{so}_{2l}$ is restricted to the subalgebra $\mathfrak{so}_{2l-1}$, the representation $S_{1/2}^{2l}$ turns into $S^{2l-1}$, and when the algebra is restricted $\mathfrak{so}_{2l+1}$ to subalgebra $\mathfrak{so}_{2l}$ space $S^{2l+1}$ turns into $S^{2l}=S_{1/2}^{2l}\oplus S_{-1/2}^{2l}$.

As will follow from the classification below, for the symplectic structure of the space $J_1$ resulting from a short $SL_2$-structure on a special Lie algebra, two different cases are possible:

1. If the representation of the Lie algebra $\mathfrak{g}_0$ on the space $J_1$ is irreducible, then according to Schur's Lemma on $J_1$ there exists a unique $\mathfrak{g}_0$-invariant non-degenerate skew-symmetric form up to a scalar multiplication. So, the symplectic structure of the space $J_1$ in this case is completely determined by the representation of $\mathfrak{g}_0$ on $J_1$. We will call this case the irreducible case of setting a symplectic structure on the space $J_1$.

2. If the representation $\mathfrak{g}_0:J_1$ is reducible, then the space $J_1$ decomposes into a direct sum of $\mathfrak{g}_0$-invariant Lagrangian subspaces, scalar multiplication between the elements of which is defined as pairing. We will call this case the Lagrangian case of setting a symplectic structure on the space $J_1$.

Let's put $\alpha_0 = -\alpha_{\operatorname{max}}$. We will mark the lowest root $\alpha_0$ of the Lie algebra $\mathfrak{g}$ on the diagram of the short $SL_2$-structure using a double circle, so as not to confuse it with zero vertices.

{\bf The case $\bf\mathfrak{g} = G_2$}. The highest root of this algebra has the form: $$\alpha_{\operatorname{max}} = 3\alpha_1 + 2\alpha_2,$$ where $\alpha_1, \alpha_2$ --- simple roots. Then the only possible case here is the one in which $\alpha_2(h) = 1, \alpha_1(h) = 0$. Consider the extended Dynkin diagram:
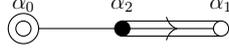
\begin{center}
\begin{dynkinDiagram} [extended, mark = {o}, labels*= {\alpha_0,\alpha_{2},\alpha_1}, root radius=.10cm, edge length=1.3cm]G{*o}
\dynkinRootMark{O}0
\end{dynkinDiagram}
\end{center}

It is clear that $\widetilde{G}^0 = SL_2$. Since $\langle\alpha_{\operatorname{max}}|\alpha_1\rangle = 0$, then $\dim\mathfrak{g}^2 =1$ and the action $\widetilde{G}^0:\mathfrak{g}^2$ is trivial, so every polynomial on $\mathfrak{g}^2$ will be invariant, that is, in this case a short $SL_2$-structure exists.

Since $\dim\mathfrak{g}^2 =1$, then $\widetilde{\mathfrak{g}}^0 = \mathfrak{g}_0 = \mathfrak{sl}_2$. The diagram shows that the representation of the algebra $\mathfrak{g}_0 = \mathfrak{sl}_2$ on the vector space $\mathfrak{g}^1$ is irreducible and the lowest weight of this representation is $\alpha_2$. From this it is easy to see that:
$$\mathfrak{g}_0 = \mathfrak{sl}_2,\quad J_1 = \operatorname{Sym}^3\mathbb{C}^2,\quad J_2 =\mathbb{C}.$$
Here the irreducible case of setting the symplectic structure of the space $J_1$ is realized.

{\bf The case  $\bf\mathfrak{g} = F_4$}.In this algebra $$\alpha_{\operatorname{max}} = 2\alpha_1 + 4\alpha_2 + 3\alpha_3 + 2\alpha_4,$$ therefore only two cases are possible:
\begin{enumerate}
\item $\alpha_1(h) = 1, \alpha_2(h) = \alpha_3(h) = \alpha_4(h) = 0$;
\item $\alpha_4(h) = 1, \alpha_1(h) = \alpha_2(h) = \alpha_3(h) = 0$.
\end{enumerate}
In case 1 , the Dynkin diagram looks like this:
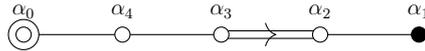
\begin{center}
\begin{dynkinDiagram} [extended, mark = {o}, labels*= {\alpha_0,\alpha_4,\alpha_3,\alpha_2,\alpha_1}, root radius=.10cm, edge length=1.3cm]F{ooo*}
\dynkinRootMark{O}0
\end{dynkinDiagram}
\end{center}

Here $\dim\mathfrak{z}^0 =1$ and $\widetilde{G}^0 =Spin_7$. The diagram shows that $\langle\alpha_{\operatorname{max}}|\alpha_4\rangle = 1, \langle\alpha_{\operatorname{max}}|\alpha_2\rangle=\langle\alpha_{\operatorname{max}}|\alpha_3\rangle = 0$, therefore, this representation $\widetilde{G}^0$ is a vector representation of the spinor group, that is, $\mathfrak{g}^2 = \mathbb{C}^{7}$. The invariant of the action of this group is obviously the scalar square of the vector, so in this case a short $SL_2$-structure exists. In this case, the Lie algebra $\tilde{\mathfrak{g}}^0$ is isomorphic to $\mathfrak{so}_7$.
The diagram also shows that the representation of the algebra $\tilde{\mathfrak{g}}^0$ on the space $\mathfrak{g}_1$ is irreducible and the lowest weight of this representation is $\alpha_1$, from which we deduce that this representation is a spinor representation of the algebra $\mathfrak{so}_7$, that is, $J_1 = S^7$.

From the reasoning described in the previous paragraphs, it follows that:
$$\mathfrak{g}^0 = \tilde{\mathfrak{g}}^0\oplus\bigl\langle\tilde{h}\otimes\mathbb{I}\bigr\rangle=\mathfrak{g}_0\oplus\bigl(\tilde{h}\otimes J_2\bigr),\quad\mathfrak{g}^2 =\tilde{e}\otimes J_2,$$
where $\mathbb{I}$ denotes the unit of the Jordan algebra $J_2$, and $\tilde{e}, \tilde{f}, \tilde{h}$ --- basic elements of the abstract Lie algebra $\mathfrak{sl}_2$ satisfying the relations $(\ref{sl2})$.

Consider the action $\tilde{\mathfrak{g}}^0:\mathfrak{g}^2$. It is clear that
$$\mathfrak{g}_0 = \bigl\{D\in\tilde{\mathfrak{g}}^0: [D,\tilde{e}\otimes\mathbb{I}] = 0\bigr\}$$ is stabilizer of the element $\tilde{e}\otimes\mathbb{I}\in\mathfrak{g}^2$ under the action of algebra $\tilde{\mathfrak{g}}^0$.
Since the Lie algebra $\mathfrak{g}_0$ is reductive, the vector $\tilde{e}\otimes\mathbb{I}$ is nonisotropic, from which it is easy to deduce that $\mathfrak{g}_0 = \mathfrak{so}_6$. When limiting the Lie algebra $\tilde{\mathfrak{g}}^0 = \mathfrak{so}_7$ to the Lie subalgebra $\mathfrak{g}_0 = \mathfrak{so}_6$, the space $S^7$ will turn into $S^6$, so $J_1 = S^6 = S^6_{1/2}\oplus S^6_{-1/2}$. The symplectic structure of the space $J_1$ corresponds to the Lagrangian case.

Similarly, when limiting the Lie algebra $\tilde{\mathfrak{g}}^0 = \mathfrak{so}_7$ to the Lie subalgebra $\mathfrak{g}_0 = \mathfrak{so}_6$, the space $J_2 =\mathbb{C}^7$ decomposes into a direct sum two orthogonal subspaces to each other: one-dimensional, spanned by the unit $\mathbb{I}$ of the algebra $J_2$, on which $\mathfrak{g}_0$ acts trivially, and six-dimensional, on which $\mathfrak{g}_0$ acts tautologically. So $J_2 = \mathbb{C}^6\oplus\mathbb{C}$.

Let's define the Jordan operation on the algebra $J_2$. Note that it is sufficient to define the Jordan multiplication between the elements of $\mathbb{C}^6$, since the additional one-dimensional subspace in $J_2$ is spanned over the unit of $J_2$. It is clear that the Jordan operation on $\mathbb{C}^6$ can be interpreted as a symmetric $\mathfrak{g}_0$-equivariant bilinear map, the image of which lies in $J_2$. Then it is passed through $$\operatorname{Sym}^2(\mathbb{C}^6) = \operatorname{Sym}^2_0(\mathbb{C}^6)\oplus\mathbb{C},$$
where the first term is an irreducible representation of the algebra $\mathfrak{so}_6$ corresponding to the doubled first fundamental weight. Then it is obvious that the result of applying the Jordan operation for vectors from $\mathbb{C}^6$ belongs to a one-dimensional subspace $J_2$ spanned over the unit of algebra $J_2$. Thus:
\begin{equation*}
A\circ B = \lambda_{A, B}\mathbb{I},\quad\forall A, B\in\mathbb{C}^6\subset J_2.
\end{equation*}
Here $\lambda_{A,B}\in\mathbb{C}$ is a proportionality coefficient uniquely determined by the pair $A, B$. Having calculated the trace of the left and right sides of the equation above, it is easy to deduce that $$\lambda_{A, B} = \dfrac{(A,B)}{\operatorname{dim} J_1}.$$
Using the formula above, the Jordan multiplication operation by $J_2$ is determined.

In all subsequent cases of short $SL_2$-structures in which $J_2\neq\mathbb{C}$ the Jordan structure will be defined in the same way, so we will not dwell on it in such detail anymore.

Summarizing the above, in this case we have:
$$\mathfrak{g}_0 = \mathfrak{so}_6,\quad J_1 = S^6_{1/2}\oplus S^6_{-1/2},\quad J_2 = \mathbb{C}^6\oplus\mathbb{C}.$$

In case 2 , the Dynkin diagram has the form:
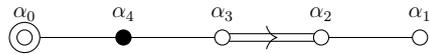
\begin{center}
\begin{dynkinDiagram} [extended, mark = {o},labels*= {\alpha_0,\alpha_4,\alpha_3,\alpha_2,\alpha_1}, root radius=.10cm, edge length=1.3cm]F{*ooo}
\dynkinRootMark{O}0
\end{dynkinDiagram}
\end{center}
In this case, it is obvious $\widetilde{G}^0 = Sp_6$, and $\dim\mathfrak{g}^2 =1$, since all Cartan numbers between the highest root and the simple roots corresponding to zero vertices are 0. Similarly to the previous reasoning, we get that in this case there is a short $SL_2$-structure, with:
\begin{equation*}
\mathfrak{g}_0 = \mathfrak{sp}_6,\quad J_1 = \sideset{}{_0^3}\bigwedge(\mathbb{C}^6),\quad J_2 = \mathbb{C}.
\end{equation*}
Here, $\sideset{}{_0^3}\bigwedge(\mathbb{C}^6)$ denotes the space of skew-symmetric tensors of the third order, whose complete contraction with the symplectic form is 0. The symplectic structure of the space $J_1$ here also corresponds to the irreducible case.

{\bf The case $\bf\mathfrak{g} = E_6$}. The highest root of this algebra has the form: $$\alpha_{\operatorname{max}} = \alpha_1 + 2\alpha_2 + 3\alpha_3 + 2\alpha_4 +\alpha_5+2\alpha_6.$$
There are 4 possible cases here:
\begin{enumerate}
\item $\alpha_2(h) = 1, \alpha_1(h) = \alpha_3(h)= ... =\alpha_6(h) = 0$;
\item $\alpha_4(h) = 1, \alpha_1(h) =...=\alpha_3(h) = \alpha_5(h) = ... =\alpha_6(h) = 0$;
\item $\alpha_6(h) = 1, \alpha_1(h) =...=\alpha_5(h) = 0$;
\item $\alpha_1(h) = \alpha_5(h) = 1, \alpha_2(h) =...=\alpha_4(h) = 0$.
\end{enumerate}
Cases 1 and 2 are also equivalent, so it is sufficient to consider only case 1. Consider it. The Dynkin diagram here has the form:
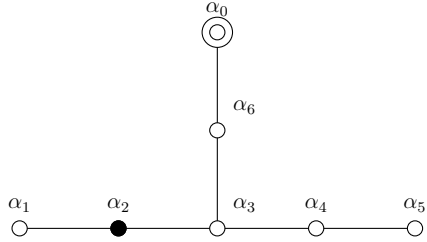
\begin{center}
\begin{dynkinDiagram}[extended, mark=o,root radius=.10cm, edge length=1.3cm,labels*= {\alpha_0,\alpha_1,\alpha_6,\alpha_2,\alpha_3,\alpha_4,\alpha_5}]E6
\dynkinRootMark{*}3
\dynkinRootMark{O}0
\end{dynkinDiagram}
\end{center}

From the diagram we get that $\dim\mathfrak{z}^0 = 1$, so $\widetilde{G}^0 = SL_2\times SL_5$. Since the highest root $\alpha_{\operatorname{max}}$ has a nonzero Cartan number only with the extreme simple root of the algebra $\mathfrak{sl}_5$ and this number is 1, then $\mathfrak{g}^2 = \mathbb{C}^5$, which says that $SL_2$ is on $\mathfrak{g}^2$ acts trivially and there are no nontrivial invariants in this case, which means there is no short $SL_2$-structure.

Let's move on to case 3 and its Dynkin diagram:
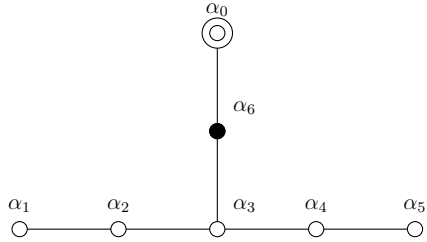
\begin{center}
\begin{dynkinDiagram}[extended, mark=o,root radius=.10cm, edge length=1.3cm, labels*= {\alpha_0,\alpha_1,\alpha_6,\alpha_2,\alpha_3,\alpha_4,\alpha_5}]E6
\dynkinRootMark{*}2
\dynkinRootMark{O}0
\end{dynkinDiagram}
\end{center}

Acting similarly to the previous cases, we get that $\widetilde{G}^0 = SL_6$ and $\dim\mathfrak{g}^2 = 1$, which means there are invariants here, that is, a short $SL_2$-structure exists, and:
\begin{equation*}
\mathfrak{g}_0 = \mathfrak{sl}_6,\quad J_1 = \sideset{}{^3}\bigwedge(\mathbb{C}^6),\quad J_2 = \mathbb{C}.
\end{equation*}
The symplectic structure of the space $J_1$ here corresponds to the irreducible case.

In the latter case , we have the following Dynkin diagram:
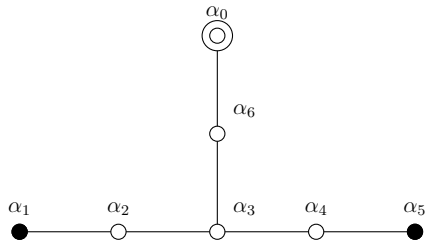
\begin{center}
\begin{dynkinDiagram}[extended, mark=o,root radius=.10cm, edge length=1.3cm, labels*= {\alpha_0,\alpha_1,\alpha_6,\alpha_2,\alpha_3,\alpha_4,\alpha_5}]E6
\dynkinRootMark{*}1
\dynkinRootMark{O}0
\dynkinRootMark{*}6
\end{dynkinDiagram}
\end{center}

Obviously, the semisimple part of $\mathfrak{g}^0$ is equal to $\mathfrak{so}_8$, but in this case $\dim\mathfrak{z}^0 =2$, that is, $\widetilde{\mathfrak{g}}^0 = \mathfrak{so}_8\oplus\langle h_1\rangle$, where $h_1\in\mathfrak{z}^0, h_1\bot h$. Consider an involutive non-identical automorphism of the Dynkin diagram, which is a symmetry with respect to the central vertical segment. This automorphism is an automorphism of the Lie algebra $\mathfrak{g}$, while it leaves the roots $\alpha_3, \alpha_6$ in place and swaps the roots $\alpha_1$ and $\alpha_5$, and also swaps $\alpha_2$ and $\alpha_4$. The elements $h$ and $h_1$ obviously define a two-dimensional space orthogonal to a four-dimensional space spanned by simple roots corresponding to zero vertices. The automorphism in question preserves this four-dimensional space, and therefore preserves the orthogonal complement. Since the marks of all vertices are preserved under the action of an automorphism, this automorphism also preserves the element $h$. So, by virtue of involution, the automorphism acts on $h_1$ by multiplying by $-1$, so we can assume that $\alpha_1(h_1) = -\alpha_5(h_1) = 1$, and all other roots on $h_1$ are equal to 0. This means that $\alpha_{\operatorname{max}}(h_1) = 0$, that is, $h_1$ acts trivially on $\mathfrak{g}^2$, so we can assume that $\widetilde{G}^0 = Spin_8$. As can be seen from the diagram $\langle\alpha_{\operatorname{max}}|\alpha_6\rangle = 1, \langle\alpha_{\operatorname{max}}|\alpha_2\rangle=\langle\alpha_{\operatorname{max}}|\alpha_3\rangle= \langle\alpha_{\operatorname{max}}|\alpha_4\rangle = 0$, that is, $\mathfrak{g}^2= \mathbb{C}^8$, which suggests that a short $SL_2$-structure exists in this case.

In this case, the representation of the algebra $\tilde{\mathfrak{g}}^0$ on the space $J_1$ is reducible, and $$J_1= S^8_{1/2}\oplus S^8_{-1/2}=S^7\oplus S^7.$$
Similarly to the previous cases, we have:
$$\mathfrak{g}_0 = \mathfrak{so}_7,\quad J_1 = S^7\oplus S^7,\quad J_2 = \mathbb{C}^7\oplus\mathbb{C}.$$
The symplectic structure of the space $J_1$ corresponds to the Lagrangian case.

{\bf The case $\bf\mathfrak{g} = E_7$}. It is known that for this algebra $$\alpha_{\operatorname{max}} = \alpha_1 +2\alpha_2+3\alpha_3 +4\alpha_4 +3\alpha_5 +2\alpha_6+2\alpha_7,$$ therefore 3 cases are possible:
\begin{enumerate}
\item $\alpha_2(h) = 1, \alpha_1(h) = \alpha_3(h) = ... = \alpha_7(h) = 0$;
\item $\alpha_6(h) = 1, \alpha_1(h) = ... = \alpha_5(h) = \alpha_7(h) = 0$;
\item $\alpha_7(h) = 1, \alpha_1(h) = ... = \alpha_6(h) = 0$.
\end{enumerate}

In case 1 , the Dynkin diagram has the form:
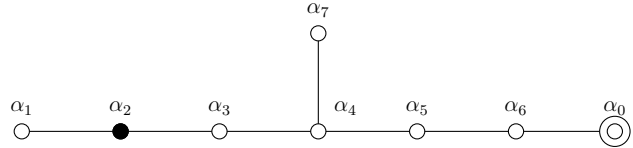
\begin{center}
\begin{dynkinDiagram}[extended, mark=o,root radius=.10cm, edge length=1.3cm, labels*= {\alpha_1,\alpha_2,\alpha_7,\alpha_3,\alpha_4,\alpha_5,\alpha_6,\alpha_0}]E7
\dynkinRootMark{*}1
\dynkinRootMark{O}7
\end{dynkinDiagram}
\end{center}

Based on it, we get that $\widetilde{G}^0 = SL_2\times Spin_{10}$. Since the highest root has a nonzero Cartan number with only the first simple root of the algebra $\mathfrak{so}_{10}$, then $\mathfrak{g}^2= \mathbb{C}^{10}$, which means that this representation has an invariant, and:
$$\mathfrak{g}_0 = \mathfrak{sl}_2\oplus\mathfrak{so}_9,\quad J_1 = \mathbb{C}^2\otimes S^9,\quad J_2 = \mathbb{C}^9\oplus\mathbb{C}.$$
The symplectic structure of the space $J_1$ corresponds to the irreducible case.

In case 2 , the Dynkin diagram has the form:
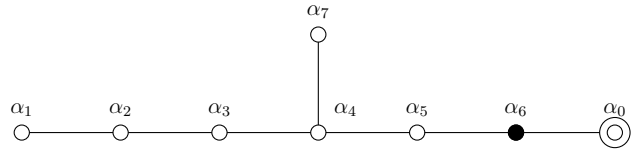
\begin{center}
\begin{dynkinDiagram}[extended, mark=o,root radius=.10cm, edge length=1.3cm, labels*= {\alpha_1,\alpha_2,\alpha_7,\alpha_3,\alpha_4,\alpha_5,\alpha_6,\alpha_0}]E7
\dynkinRootMark{*}6
\dynkinRootMark{O}7
\end{dynkinDiagram}
\end{center}

Here $\widetilde{G}^0= Spin_{12}$, and $\dim\mathfrak{g}^2=1$, so there is an invariant and
$$\mathfrak{g}_0 = \mathfrak{so}_{12},\quad J_1 = S^{12}_{1/2},\quad J_2 = \mathbb{C}.$$
The symplectic structure of the space $J_1$ here corresponds to the irreducible case.

In case 3 , the Dynkin diagram has the form:
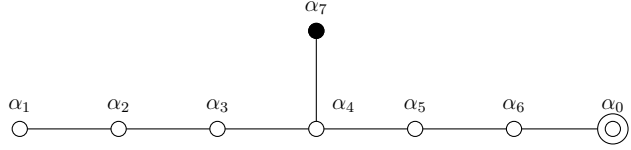
\begin{center}
\begin{dynkinDiagram}[extended, mark=o,root radius=.10cm, edge length=1.3cm, labels*= {\alpha_1,\alpha_2,\alpha_7,\alpha_3,\alpha_4,\alpha_5,\alpha_6,\alpha_0}]E7
\dynkinRootMark{*}2
\dynkinRootMark{O}7
\end{dynkinDiagram}
\end{center}

Here $\widetilde{G}^0 = SL_7$, and $\mathfrak{g}^2 = \mathbb{C}^{7}$, so there are no nontrivial invariants for this action.

{\bf The case $\bf\mathfrak{g} = E_8$}. In this algebra, the highest root has the form: $$\alpha_{\operatorname{max}} = 2\alpha_1 + 3\alpha_2 + 4\alpha_3 + 5\alpha_4 + 6\alpha_5 + 4\alpha_6 + 2\alpha_7 + 3\alpha_8,$$ therefore, only 2 cases are possible here:
\begin{enumerate}
\item $\alpha_1(h) = 1, \alpha_2(h) = ... = \alpha_8(h) = 0$;
\item $\alpha_7(h) = 1, \alpha_1(h) = ... = \alpha_6(h) = \alpha_7(h) = 0$;
\end{enumerate}
The first case has the following Dynkin diagram:
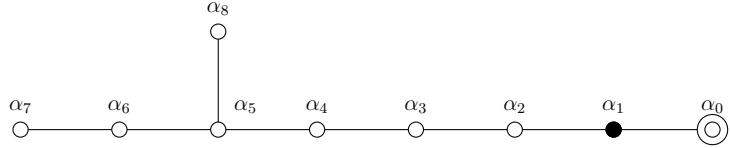
\begin{center}
\begin{dynkinDiagram}[extended, mark=o,root radius=.10cm, edge length=1.3cm, labels*= {\alpha_0,\alpha_7,\alpha_8,\alpha_6,\alpha_5,\alpha_4,\alpha_3,\alpha_2,\alpha_1}]E8
\dynkinRootMark{*}8
\dynkinRootMark{O}0
\end{dynkinDiagram}
\end{center}
In this case, $\widetilde{G}^0 = E_7$, and $\dim\mathfrak{g}^2 =1$, that is, this representation has an invariant and
$$\mathfrak{g}_0 = E_7,\quad J_1 = V(\pi_1),\quad J_2 = \mathbb{C}.$$
Here $V(\pi_1)$ denotes the representation space of the first fundamental weight of the algebra $E_7$. The symplectic structure of the space $J_1$ corresponds to the irreducible case.
\newline
\newline
The Dynkin diagram of the second case has the form:
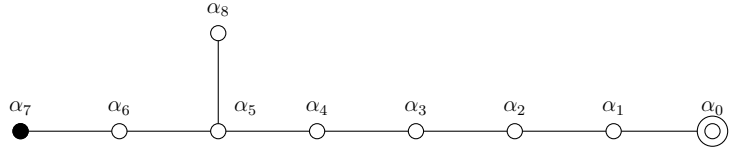
\begin{center}
\begin{dynkinDiagram}[extended, mark=o,root radius=.10cm, edge length=1.3cm, labels*= {\alpha_0,\alpha_7,\alpha_8,\alpha_6,\alpha_5,\alpha_4,\alpha_3,\alpha_2,\alpha_1}]E8
\dynkinRootMark{*}1
\dynkinRootMark{O}0
\end{dynkinDiagram}
\end{center}

Similarly to the previous arguments, we get that $\widetilde{G}^0 = Spin_{14}$, and $\mathfrak{g}^2$ is its vector representation, that is, this representation has an invariant, moreover:
$$\mathfrak{g}_0 = \mathfrak{so}_{13},\quad J_1=S^{13},\quad J_2 =\mathbb{C}^{13}\oplus\mathbb{C}.$$
The symplectic structure of the space $J_1$ corresponds to the irreducible case.

Thus, all existing short $SL_2$-structures can be listed in the following table:
\begin{center}
\scalebox{0.55}{
\begin{tabular}{ ||c | c | c | c | c || }
\hline
$\mathfrak{g}$&Dynkin 's scheme&$J_1$&$\mathfrak{g}_0$ & $J_2$ \\ \hline
$\mathfrak{sl}_n$&\begin{dynkinDiagram} [text style={scale=0.8,blue},labels*= {,,,\alpha_{i},,,\alpha_{j}},root radius=.10cm, edge length=1.3cm]A{oo.o*o.o*o.oo}
\end{dynkinDiagram}, $i = n-j$
 &$\begin{array}{c}\left(\mathbb{C}^{i}\otimes\bigl(\mathbb{C}^{n-2i}\bigr)^{*}\right)\oplus\\
 \oplus\left(\mathbb{C}^{n-2i}\otimes\bigl(\mathbb{C}^{i}\right)^{*}\bigr)\end{array}$& $\mathfrak{sl}_i\oplus\mathfrak{sl}_{n-2i}\oplus\mathfrak{c}_{i,n}$ & $\mathfrak{gl}_i$ \\ \hline
$\mathfrak{so}_{2n+1}$&\begin{dynkinDiagram} [text style={scale=0.8,blue},labels*= {,,,\alpha_{i}},root radius=.10cm, edge length=1.3cm]B{oo.o*o.oo}
\end{dynkinDiagram}, $i=2k$ &$\mathbb{C}^{i}\otimes\bigl(\mathbb{C}^{2(n-i)+1}\bigr)^{*}$& $\mathfrak{sp}_i\oplus\mathfrak{so}_{2(n-i)+1}$ &$\mathfrak{sym}^{-}_i$ \\ \hline
$\mathfrak{so}_{2n}$&\begin{dynkinDiagram} [text style={scale=0.8,blue},labels*= {,,,\alpha_{i}},root radius=.10cm, edge length=1.3cm]D{oo.o*o.ooo}
 \end{dynkinDiagram}, $i=2k\leqslant n-2$&$\mathbb{C}^{i}\otimes\bigl(\mathbb{C}^{2(n-i)}\bigr)^{*}$&  $\mathfrak{sp}_i\oplus\mathfrak{so}_{2(n-i)}$ & $\mathfrak{sym}^{-}_i$ \\ \hline
$\mathfrak{so}_{2n}$&\begin{dynkinDiagram} [text style={scale=0.8,blue},labels*= {,,,\alpha_{n-1},\alpha_n},root radius=.10cm, edge length=1.3cm]D{oo.o**}
\end{dynkinDiagram}, $n=2k+1$&$\mathbb{C}^{n-1}\otimes\bigl(\mathbb{C}^{2}\bigr)^{*}$&$\mathfrak{sp}_{n-1}\oplus\mathfrak{so}_2$& $\mathfrak{sym}^{-}_{n-1}$\\ \hline
$\mathfrak{sp}_{2n}$&
\begin{dynkinDiagram} [text style={scale=0.8,blue},labels*= {,,,\alpha_{i}},root radius=.10cm, edge length=1.3cm]C{oo.o*o.oo}
\end{dynkinDiagram}, $i<n$ &$\mathbb{C}^{i}\otimes\bigl(\mathbb{C}^{2(n-i)}\bigr)^{*}$& $\mathfrak{so}_{i}\oplus \mathfrak{sp}_{2(n-i)}$&$\mathfrak{sym}^{+}_{i}$\\ \hline
$G_2$&\begin{dynkinDiagram} [text style={scale=0.8,blue}, mark = {o},labels*= {\alpha_{2},\alpha_1}, root radius=.10cm, edge length=1.3cm]G{*o}
\end{dynkinDiagram}&$\operatorname{Sym}^3\mathbb{C}^2$&$\mathfrak{sl}_2$&$\mathbb{C}$
  \\ \hline
$F_4$&\begin{dynkinDiagram} [text style={scale=0.8,blue},mark = {o},labels*= {\alpha_4,\alpha_3,\alpha_2,\alpha_1}, root radius=.10cm, edge length=1.3cm]F{ooo*}
\end{dynkinDiagram}&$S^6_{1/2}\oplus S^6_{-1/2}$&$\mathfrak{so}_6$&$\mathbb{C}^6\oplus\mathbb{C}$\\ \hline
$F_4$&\begin{dynkinDiagram} [text style={scale=0.8,blue}, mark = {o},labels*= {\alpha_4,\alpha_3,\alpha_2,\alpha_1}, root radius=.10cm, edge length=1.3cm]F{*ooo}
\end{dynkinDiagram}&$\bigwedge_0^{3}\left(\mathbb{C}^6\right)$&$\mathfrak{sp}_6$&$\mathbb{C}$\\ \hline
$E_6$&\begin{dynkinDiagram}[text style={scale=0.8,blue}, mark={o},root radius=.10cm, edge length=1.3cm,labels*= {\alpha_1,\alpha_6,\alpha_2,\alpha_3,\alpha_4,\alpha_5}]E6
\dynkinRootMark{*}2
\end{dynkinDiagram}&$\bigwedge^{3}\left(\mathbb{C}^6\right)$&$\mathfrak{sl}_6$&$\mathbb{C}$\\ \hline
$E_6$&\begin{dynkinDiagram}[text style={scale=0.8,blue}, mark={o},root radius=.10cm, edge length=1.3cm,labels*= {\alpha_1,\alpha_6,\alpha_2,\alpha_3,\alpha_4,\alpha_5}]E6
\dynkinRootMark{*}1
\dynkinRootMark{*}6
\end{dynkinDiagram}&$S^{7}\oplus S^{7}$&$\mathfrak{so}_7$&$\mathbb{C}^7\oplus\mathbb{C}$\\ \hline
$E_7$&\begin{dynkinDiagram}[text style={scale=0.8,blue},mark={o},root radius=.10cm, edge length=1.3cm,labels*= {\alpha_6,\alpha_7,\alpha_5,\alpha_4,\alpha_3,\alpha_2,\alpha_1}]E7
\dynkinRootMark{*}6
\end{dynkinDiagram}&$\mathbb{C}^2\otimes S^{9}$&$\mathfrak{sl}_2\oplus\mathfrak{so}_9$&$\mathbb{C}^9\oplus\mathbb{C}$\\ \hline
$E_7$&\begin{dynkinDiagram}[text style={scale=0.8,blue},mark=o,root radius=.10cm, edge length=1.3cm,labels*= {\alpha_6,\alpha_7,\alpha_5,\alpha_4,\alpha_3,\alpha_2,\alpha_1}]E7
\dynkinRootMark{*}1
\end{dynkinDiagram}&$S_{1/2}^{12}$&$\mathfrak{so}_{12}$&$\mathbb{C}$\\ \hline
$E_8$&\begin{dynkinDiagram}[text style={scale=0.8,blue}, mark=o,root radius=.10cm, edge length=1.3cm,labels*= {\alpha_7,\alpha_8,\alpha_6,\alpha_5,\alpha_4,\alpha_3,\alpha_2,\alpha_1}]E8
\dynkinRootMark{*}8
\end{dynkinDiagram}&$V(\pi_1)$&$E_7$&$\mathbb{C}$\\ \hline
$E_8$&\begin{dynkinDiagram}[text style={scale=0.8,blue}, mark=o,root radius=.10cm, edge length=1.3cm,labels*= {\alpha_7,\alpha_8,\alpha_6,\alpha_5,\alpha_4,\alpha_3,\alpha_2,\alpha_1}]E8
\dynkinRootMark{*}1
\end{dynkinDiagram}&$S^{13}$&$\mathfrak{so}_{13}$&$\mathbb{C}^{13}\oplus\mathbb{C}$\\
\hline
\end{tabular}
}
\end{center}
Note that the Lie algebra $\mathfrak{g}_0$ is not completely defined by the pair $(J_1, J_2)$. For example, in the case of a simple Lie algebra $\mathfrak{g}=E_8$ for $J_1 = V(\pi_1)$. The Jordan algebra $J_2$ is one-dimensional, as in the case of the algebra $\mathfrak{g} = \mathfrak{sl}_n$ for $i=1$. If we match $n$ so that $\dim V(\pi_1) = \dim J_1 (\mathfrak{sl}_n)$, where $J_1 (\mathfrak{sl}_n) = \left(\mathbb{C}^{i}\otimes\left(\mathbb{C}^{n-2i}\right)^{*}\right)\oplus\left(\mathbb{C}^{n-2i}\otimes\left(\mathbb{C}^{i}\right)^{*}\right)$, then it can be seen that the algebras $\mathfrak{g}_0$ are different in these cases, while the pairs $(J_1, J_2)$ are coincide. Thus, for the pair $(J_1, J_2)$ it is impossible to uniquely construct a simple symplectic Lie-Jordan structure.
\bibliographystyle{amsplain}

\end{document}